\documentstyle[12pt]{article}    
\newcounter{props}[section]
\newcommand{\wrt}{with respect to\ }

\newcommand{\tx}{there exist\ }
\newcommand{\txs}{there exists\ }
\newcommand{\st}{such that\ }
\newcommand{\fe}{ for every\  }

\newcommand{\sq}{sequence}
\newcommand{\lra}{\longrightarrow }

\newcommand{\ty}{\infty }

\newcommand{\calb}{{\cal B}}

\newcommand{\calg}{\mbox{$\cal{G}$}}

\newcommand{\calt}{\mbox{$\cal{T}$}}

\newcommand{\nin}{n\in{\bf N}}

\newcommand{\NN} {\mbox{${\bf N}$}}\newcommand{\nn}{{\bf N}}
\newcommand{\N} {\mbox{$[{\bf N}]$}}

\newcommand{\bigcupnf}{\bigcup_{n=1}^{\infty}}
\newcommand{\sumnf}{\sum_{n=1}^{\infty}}
\newcommand{\sumin}{\sum_{i=1}^{n}}

\newcounter{excounter}

\begin{document}           
\newcommand{\sprt}{\vspace{0.15in}}
\newcommand{\sprh}{\vspace{0.10in}}
\newcommand{\beq}{\begin{eqnarray*}}
\newcommand{\eeq}{\end{eqnarray*}}
\newcommand{\bce}{\centerline}
\newcommand{\ece}{ }
\newcommand{\hardexercise}{
\hspace{-2.5mm}$ ^{*}$\hspace{1.5mm}}
\newcommand{\spp}{{\rm supp}\,}
\newcommand{\sbs}{\subset}
\newcommand{\sps}{\supset}
\newcommand{\stm}{\setminus}
\newcommand{\lf}{\left}\newcommand{\r}{\right}
\newcommand{\f}{\frac}
\newcommand{\DD}{\displaystyle}
\newcommand{\TT}{\textstyle}
\newcommand{\th}[1]{\sprh{\bf #1:}}
\newcommand{\proof}{
  \sprh \bf Proof.  \rm \  }
\newcommand{\thm}[1]{\stepcounter{props}
  \sprt\bf \noindent \thesection .\theprops .\ 
  #1:\rm \  }
\newcommand{\sol}[1]{
  \sprt\hspace{1cm}\bf  #1\par\hspace{-5pt}\nopagebreak\rm \  }
\newcommand{\exercises}{\sprt
\[ {\rm EXERCISES}\]\sprt\small
\setcounter{excounter}{0}\begin{list}%
{\thesection --\theexcounter}{\usecounter{excounter}
}}
\newcommand{\de}{\delta}
\newcommand{\e}{\epsilon}\newcommand{\vare}{\varepsilon}
\newcommand{\om}{\omega}\newcommand{\z}{\zeta}
\newcommand{\abst}{\sprt
\[ {\rm ABSTRACT}\]\sprt
\begin{description}
\rightmargin=20 mm}
\newcommand{\stareq}[3]{\vspace{7pt}$(#1)\hspace{#2cm}{\displaystyle
#3}$ \vspace{7pt}\\}
\newcommand{\bb}{
                 \item [ }
\newcommand{\ce}[2]{\vspace{7pt}\noindent \mbox{}\hfill
${\displaystyle #1}$ \hfill #2\vspace{7pt}\noindent}
\frenchspacing
\newcommand{\calf}{{\mbox{ \boldmath$\cal F$}}}
\newcommand{\cf}  {{\mbox{ \boldmath$\cal F$}}}
\newcommand{\xii}{{\mbox{ \boldmath$\xi$}}}
\newcommand{\zi}{{\mbox{ \boldmath$\z$}}}
\newcommand{\zii}{{\mbox{\boldmath$\z$}}}
\newcommand{\zin}{{\mbox{\boldmath$\zi+n$}}}
\newcommand{\zinn}{{\mbox{\boldmath$\zi+n-1$}}}
\newcommand{\etai}{{\mbox{ \boldmath$\eta$}}}

\Large
\begin{center}\bf Convex unconditionality and summability
of weakly null  sequences
\end{center}

\normalsize
\[{\rm by}\]
\[ {\rm S.\; A.\; Argyros,\;\; S.\ Merkourakis\;\; and\;\; A.\; Tsarpalias}\]
\[ {\rm (Athens,\;\; Greece)}\]

\sprt\sprt
{\sc Abstract.}
{\em It is proved that every normalized weakly null \sq\ has a sub\sq\
which is convexly unconditional. Further, an Hierarchy of
summability methods is introduced and with this we give a complete
classification of the complexity of weakly null \sq s.}

\frenchspacing
\sprt
\sprt
{\bf INTRODUCTION}
\normalsize
In the present paper we investigate the behavior of the sub\sq s
of a weakly null \sq\ $(x_n)_{\nin}$ of a Banach space $X$ \wrt   two
fundamental properties. The first is
the convex unconditionality which is investigated in the
first section of the paper. This is defined as:

\th{Definition:} A normalized \sq\ $(x_n)_{n\in{\bf N}}$ in a Banach
space $X$ is said to be {\em convexly unconditional\/} if \fe\ 
$\delta>0$ \txs $C(\de)>0$ \st if an absolutely convex combination
$x=\sumnf a_nx_n$ satisfies $\|x\|>\de$ then
$\|\sumnf \varepsilon_na_nx_n\|>C(\de)$ \fe
choice of signs $(\varepsilon_n)_{n\in{\bf N}}$.

\sprh
The result we prove here is the following theorem.

\th{Theorem A} If $(x_n)_{\nin}$ is a normalized weakly null \sq\
 in a Banach space $X$ then it has a convexly unconditional
subsequence.

\sprh
A fundamental example due to B. Maurey and H. Rosenthal [M--R] showed
that we could not expect that every normalized weakly null \sq\
has an unconditional sub\sq. The recent examples [G--M], [A--D]
show that there are spaces without any unconditional basic \sq.
On the other hand there are results where some weaker forms of
unconditionality appear. One of them is due to J. Elton [E], [O$_1$]
which is related to the unconditional behavior of the linear
combinations with coefficients away from zero 
and the other is due to E. Odell
and it is related to the unconditionality of Schreier admissible
linear combinations.
Our theorem is in the same direction with Elton's Theorem, more
precisely it is the dual result, and the proof is based, as his
proof, on infinite Ramsey Theorem.
The result follows from a combinatorial principle (Lemma 1.2) which
seems of independent interest and it is also used in the second
part of the paper.

The existence of a convexly unconditional \sq\ is a strong evidence
that the convex sets behave much better, \wrt the unconditionality
than the subspaces of a Banach space.

In the second part we deal with summability methods. The starting point
for our investigation is the following question.

As it follows from Mazur's theorem,   every weakly null \sq\ has
convex combinations norm converging to zero. The general question
is to describe ``regular'' convex combinations with this property.
This problem was faced from the early days of the development of
Banach space theory. Thus Banach and Saks proved that every $(x_n)_{\nin}$
bounded \sq\  in ${\rm L}^p(\mu)$, $1<p<\ty$ has a norm
Cesaro summable sub\sq. This result was extended by W. Szlenk 
for weakly convergent \sq s in ${\rm L}^1(\mu)$. Shortly after
Banach-Saks Theorem, an example was given by J. Schreier [Sch]
of  a weakly
null \sq\ with no norm Cesaro summable sub\sq.
Schreier's example is defined as follows:
First we define the following family.
\[\calf=\{ F\in\nn:\# F\leq\min F\}.\] 
Then on the vector space $c_{00}(\nn)$ of eventually zero \sq s
we define the norm
\[ \lf\| (a_n)_{\nin}\r\|=\sup\lf\{\sum_{n\in F} |a_n|:
F\in\calf\r\}.\]
It is easy to see that $\calf$ is compact in the topology of
pointwise convergence.
Hence the standard basis $(e_n)_{\nin}$ is weakly null. Further
from the definition of $\cf$ we get that \fe
\[n_1<n_2<\cdots<n_k\;{\rm we\; have}\;\;
\lf\|\f{e_{n_1}+e_{n_2}+\cdots+e_{n_k}}{k}\r\|\geq\f{1}{2}.\]

\noindent So no sub\sq\ of $(e_n)_{\nin}$ is norm Cesaro summable.
Later it is proved by H. Rosenthal 
that if $(x_n)_{\nin}$ is weakly null and no
sub\sq\ is norm Cesaro summable then \txs $(n_i)_{i\in \nn}$ and
$\e>0$ \st
\[\|\sum_{i\in F}a_ie_{n_i}\|>\e\cdot\sum_{i\in F}|a_i|\]
for all $F\in\calf$. 
Whenever this property appears, we say that the \sq\
$(x_{n_i})_{i\in\nn}$ is  an $\ell^1$ spreading model.
This result, in connection with a theorem proved by P. Erd\"os and
M. Magidor [E-M] gives the  following dichotomy.

\th{Theorem} For every
 $(x_n)_{\nin}$ weakly null \sq\ exactly one of the following
holds:

(a) For every $M\in\N$ \txs $L\in[M]$ \st for all $P\in[L]$
$P=(n_i)_{i\in\nn}$, the sub\sq\
 $(x_{n_i})_{i\in\nn}$ is norm Cesaro summable

(b) There exists $M\in\N$ \st $M=(m_i)_{i\in \nn}$ and the sub\sq\
$(x_{m_i})_{i\in\nn}$ is an $\ell^1$ spreading model.

\sprh 
A proof of this theorem is also given in [M].

This theorem is sufficient when condition (a) appears.
If (b) holds then there is no information on the structure of
convex combinations that converge in norm to zero. Our
aim is to give a
 full extension of the above theorem
and through this to describe the complexity of weakly null \sq s.
For this we use two hierarchies, {\em Schreier Hierarchy}
and the {\em Repeated Averages Hierarchy.}

\sprh{\bf Schreier Hierarchy:}
Schreier family $\calf$ is quite important in the theory of Banach
space. Recall that it is one of the main ingredients in the definition of
Tsirelson's space [T].  D. Alspach and S. Argyros [A-A] defined
a family $\{\calf_\xi\}_{\xi<\om_1}$ called generalized
Schreier families. The definition of $\calf_\xi$ is given
in the following way:

Set $\calf_0=\{\{n\}:\nin\}$ and $\cf_1=\cf$.

If $\cf_\xi$ has been defined then we set

\bce{$\cf_{\xi+1}=\lf\{\bigcup_{i=1}^nF_i:n\leq F_1<F_2<\cdots< F_n,
\: F_i\in\cf_\xi\r\}$.}

If $\xi$ is a limit ordinal choose $(\xi_n)_{\nin}$ strictly
increasing to $\xi$ and we set

\bce{$\cf_\xi=\{F:F\in\cf_{\xi_n}, n\leq \min F\}.$}

We decided to call this family Schreier Hierarchy since it carries
some strong universal properties some of which are described
in the present paper. Roughly speaking, the complexity of every
compact countable metric space is dominated by 
some member of $\{\cf_\xi\}_{\xi<\om_1}$.
Further members of $\{\cf_\xi\}_{\xi<\om_1}$ appear naturally
in several cases. For example, the $n^{\rm th}$ norm in the
inductive definition of Tsirelson's space is implicitly connected
to the family $\cf_n$. Explicitly the family $(\cf_n)_{\nin}$ appeared
for the first time in an example constructed by E. Odell [A-O].
Recently $\{\cf_\xi\}_{\xi<\om_1}$ are used in the investigation of
 asymptotic $\ell^p$ spaces. 
Connected to the family $\{\cf_\xi\}_{\xi<\om_1}$ is the
following definition.

\th{Definition} Let $(x_n)_{\nin}$ be a bounded \sq\ in a Banach
space $X$. For $M\in\N$, 
$M=(m_i)_{i\in\nn}$ we say that $(x_{m_i})_{i\in\nn}$
is an $\ell^1_\xi$ spreading model if \txs $\e>0$ \st for all
choices,  $(a_i)_{i\in F}$, for $F\in\cf_\xi$ we have that:

\bce{$\|\sum_{i\in F}a_ix_{m_i}\|\geq\e\cdot\sum_{i\in F}|a_i|$.}

\sprh
It is clear that $\ell^1_1$ spreading model is the usual
$\ell^1$ spreading model.
Since the families $(\cf_\xi)_{\xi<\om_1}$ are of increasing complexity,
the existence of a sub\sq\ which is an $\ell^1_\xi$ model, for
large $\xi$, describes strong $\ell^1$ behavior of the given
\sq s. As it is proved in [A-A], if a \sq\ contains $\ell^1_\xi$ spreading
models for all $\xi<\om_1$ then actually contains a sub\sq\ 
equivalent to the usual basis of $\ell^1$.

The second hierarchy introduced here is that of Repeated Averages.

{\bf Repeated Averages Hierarchy:} To introduce this we give some notations 
and definitions.

We denote by $S^+_{\ell^1}$ the positive part of the unit 
sphere of $\ell^1(\nn)$ and if $H=(x_n)_{\nin}$ is a bounded
\sq\ in a Banach space, $A=(a_n)_{\nin}\in S^+_{\ell^1}$
we set

\bce{ $A\cdot H=\sum_{n=1}^\ty a_nx_n\in X$.}

For $M\in[\nn]$ a \sq\ $(A_n)_{\nin}$ of successive blocks
in $S^+_{\ell^1}$ defines an $M$-summability method, denoted by
$M{\rm -}(A_n)_{\nin}$ if 
$M=\bigcupnf\spp A_n$.

\sprh
{\bf Definition:} A weakly null \sq\ $H=(x_n)_{\nin}$ is 
$M{\rm -}(A_n)_{\nin}$ sum\-mable if the \sq\ $(A_n\cdot H)_{\nin}$ is
Cesaro summable.

The RA Hierarchy is defined, inductively, 
\fe $M\in\N$ and $\xi<\om_1$ and it is an $M$-summability method
denoted by $(\xii^M_n)_{\nin}$.
We also use the notation  $(M,\xi)$ for the same method.
Thus the RA Hierarchy is the family
\[ \{ (M,\xi):M\in\N,\; \xi<\om_1\}.\]
The precise definition is given at the beginning of the second section
of the paper. A brief desciption of it goes as follows:
For $\xi=0$ and $M=(m_n)_{\nin}$ we set
$\xii^M_n=e_{m_n}$. Thus the $(M,\xi)$-summability, for $\xi=0$, of
a weakly null \sq\ $(x_n)_{\nin}$ is exactly the norm Cesaro
summability of the sub\sq\ $(x_{m_n})_{\nin}$ where
$M=(m_n)_{\nin}$

If $(\xii^M_n)_{\nin}$ has been defined then for $\z=\xi+1$ we set
$\zi^M_n$ to be the average of an appropriate number of
successive elements of 
$(\xii^M_n)_{\nin}$. This justifies the term Repeated Averages.
For $\z$ limit ordinal  $(\zi^M_n)_{\nin}$ is constructed by
a careful choice of terms of $\lf\{(\xii^M_n):\xi<\z, {\nin}\r\}$.

One property we would like to mention here is that
$\spp\xii^M_n\in\cf_\xi$ and moreover it is a maximal element
of $\cf_\xi$. Thus $(\xii^M_n)_{\nin}$ exhausts the complexity
of the family  $\cf_\xi$. More important is that
$(\xii^M_n)_{\nin}$ carries some nice stability properties
(see P.3 -- P.4 after the Definition) which allows us to handle them in
the proofs of the theorems.

The difference between the RA Hierarchy and the summability methods
described as an infinite matrix is that in RA Hierarchy the
summability of a sub\sq\ $(x_n)_{n\in M}$ depends on the subset
$M$ while in the usual case, after reordering $(x_n)_{n\in M}$
as $(x_{n_k})_{k\in\nn}$, we ignore the set $M$ and apply the
summability method \wrt the index $k$.
Thus in our case for fixed countable ordinal $\xi$ we have $2^\om$
methods $\{(M,\xi):M\in[\nn]\}$ which have uniformly bounded
complexity. This is so, since \fe $M\in[\nn]$, $\nin$, the set
$\spp\xii^M_n$ belongs to the compact family $\cf_\xi$.

For a given $M\in[\nn]$ the methods $\{(M,\xi):\xi<\om_1\}$ are
increasing very fast. It is worthwhile to remark that if
for $\xi<\om_1$ and $\nin$ we set
$k^\xi_n=\max\spp\xi^{\nn}_n$ then the family
$\{(k^\xi_n):\nin,\: \xi<\om_1\}$ is the Ackerman Hierarchy, a
well known hierarchy of Mathematical Logic.

\th{Theorem B} For $(x_n)_{\nin}$ weakly null \sq\ in a Banach space
$X$ and $\xi<\om_1$  exactly one of the following holds.

(a) For every $M\in[\nn]$ \txs $L\in[M]$ \st
\fe $P\in[L]$ the  \sq\ $(x_n)_{\nin}$ is $(P,\xi)$ summable

(b) There exists $M\in\N$ $M=(m_i)_{i\in\nn}$ \st
$(x_{m_i})_{i\in\nn}$ is an $\ell^1_{\xi+1}$ spreading model.

\sprh It is proved in [A-A] that every $(x_n)_{\nin}$  weakly null
\sq\ \txs $\xi<\om_1$ \st \fe $\z\geq\xi$ no sub\sq\ of
$(x_n)_{\nin}$ is an $\ell^1_\z$ spreading model. 
So we introduce the {\em Banach-Saks} index of a weakly null \sq\
defined as

$BS[(x_n)_{\nin}]=\min\{\xi:$ no sub\sq\ of $(x_n)_{\nin}$ is
an $\ell^1_\xi$ spreading model$\}$

\noindent and from Theorem B we get the following

\th{Theorem C} Let $H=(x_n)_{\nin}$ be a weakly null \sq\ with\\
$BS[(x_n)_{\nin}]=\xi$.
Then: $\xi$ is the unique ordinal satisfying the following

(a) For every $M\in[\nn]$ \txs $L\in[M]$ \st \fe
$P\in[L]$ $\lim_{\nin}\|\xii^P_n\cdot H\|=0$.

(b) For every $\z<\xi$ \txs $L_\z\in\N$ \st $L_\z=(n_i)_{i\in\nn}$ and
$(x_{n_i})_{i\in\nn}$ is an $\ell^1_\z$ spreading model.

(c) If $\xi=\z+1$ \txs $\e>0$ and $L\in\N$ \st for all
$P\in[L]$ $\lf\|\zi^P_n\cdot H\r\|>\e$ and
$(\lf(\zi^P_n\cdot H\r)_{\nin}$ is Cesaro summable.

\sprh 
For $\xi=0$ Theorem B implies exactly the dichotomy mentioned
at the beginning of the introduction (Theorem). Theorem C
gives the full description of the norm summability
for a weakly null \sq\ in terms of the methods $\{(M,\xi):
M\in[\nn],\xi<\om_1\}$. 
This justifies the universal character of these summability
methods as well as the universal character of Schreier Hierarchy since,
as we mentioned above, the $\spp\xii^M_n$ belongs to $\cf_\xi$.

\th{Definition} (a) A Banach space $X$ has the {\em$\xi$-Banach Saks property}
($\xi$-BS) if \fe bounded \sq\ $(x_n)_{\nin}$ in $X$
\txs $L\in[\nn]$ \st $(x_n)_{\nin}$ is $(L,\xi)$ summable.

(b) The space $X$ has the {\em weak $\xi$ Banach Saks property}
(w $\xi$-BS) if the above property holds only for
weakly convergent \sq s.

\sprh From Theorem B and C follow the next corollaries.

\th{Corollary} For every separable reflexive Banach space $X$
\txs a unique ordinal $\xi<\om_1$ \st

(i) For all ordinals $\z\geq \xi$ the space $X$ has $\z$-BS.

(ii) For every $\z<\xi$ the space $X$ fails $\z$-BS.

\th{Corollary} If $X$ is a separable Banach space not containing
isomorphically $\ell^1$ then \tx a unique ordinal 
$\xi<\om_1$ \st

(i) For all ordinals $\z\geq \xi$ the space $X$ has w $\z$-BS.

(ii) For every $\z<\xi$ the space $X$ fails w $\z$-BS.

\sprh
The proofs of the above theorems use infinite Ramsey theorem
and an index introduced here for compact families of
infinite subsets of $\nn$ that is called strong Cantor Bendixson
index. This index helps us to develop a criterion for embedding 
 the family $\cf_\xi$ into a family $\cf$ provided the index 
of $\cf$ is greater than $\om^\xi$.
Also, the proofs of these theorems make use of Lemma 1.2
and a variation of it.

\sprh{\bf Notation:} For $N$ infinite subset of \NN\ we denote by
$[N]$ the set of all infinite subsets of $N$.
Further, we denote by $[N]^{<\om}$ the set of all finite subsets
of the set $N$. In the sequel for $F\in[\nn]^{<\om}$ we will
identify the set $F$ with its characteristic function.
Thus for $A=(a_n)_{\nin}$ in $\ell^1(\nn)$ we will denote by
$\langle A,F\rangle$ the quantity $\sum_{n\in F} a_n$.
For $M\in[\nn]$ we will denote by $M=(m_i)_{i\in\nn}$ the
natural order of the set $M$.

\sprh 
As we mentioned above, our proofs use in an essential way the infinite
Ramsey Theorem. This theorem, one of the most important principles
in infinite combinatorics proved in several steps by Nash-Williams
[N-W], Galvin and Prikry [G-P] and in the final form by Silver
[Si]. Silver's proof was model-theoretic. Later Ellenduck [Ell]  gave
a  proof of Silver's result using classical methods. We recall the statement
of the theorem.

\th{0.1. Theorem} 
Let $A$ be an analytic subset of $[\nn]$. Then \fe $M\in[\nn]$
 \txs $L\in[M]$ \st either
$[L]\sbs A$ or else $[L]\sbs [M]\stm A$.

\sprh 
In the sequel any set $A$ satisfying the above property will be
called completely Ramsey.
Here we consider the elements of $[M]$ as strictly increasing \sq s
and we topologize $[M]$ by the topology of the pointwise convergence.

\sprt\bce{{\large \bf 1. Convex unconditionality}}

\sprh
We start with two combinatorial lemmas.

\th{1.1. Lemma} Let $F$ be a relatively weakly compact subset of $c_0(\NN)$.
Then \fe $N'\in[\NN]$ \txs $M\in[N']$ such that:
If $l_1<l_2<\cdots <l_n$ are elements of $M$ and $f\in F$ \st
\fe $i=2,\ldots,n$ $f(l_i)\geq\de$ then \txs $g\in F$ \st \fe $i=2,
\ldots, n$ $g(l_i)\geq\de$ and $|g(l_1)|<\e$.

\proof For $\nin$ we set

\vspace{5pt}
\noindent$S_n=\{$\begin{minipage}[t]{12cm}$M\in [N']:M=(m_i)$ and if
\txs $f\in F$ \st

 $\forall\; i=2,\ldots, n$ $f(m_i)\geq\de$
then \txs $g\in F$ \st

 $\forall\; i=2,\ldots, n$ $g(m_i)\geq\de$
and $g(m_1)<\e\}$\end{minipage}

\vspace{5pt}
\noindent It is clear that each $S_n$ is closed in the topology of
pointwise convergence. Hence $S=\bigcap_{n=1}^\ty S_n$ is closed and
therefore completely Ramsey. Choose $M\in[N']$ \st
either $M\sbs S$ or $[M]\sbs[N']\stm S$.

Suppose that $[M]\sbs[N']\stm S$. Then $M=(m_i)$ and consider any $\nin$.
For $1\leq j\leq n$ set
\[ L_j=\{ m_j, m_{n+1}, m_{n+2},\ldots\},\]
which does not belong to $S$. Therefore, for any such $j$ there are
$f_j\in F$ and $l_j\in\NN$ \st for all $i=1,\ldots,l_j$ $f_j(m_{n+i})
\geq\de$ and every $g\in F$ with $g(m_{n+i})\geq\de$
satisfies $|g(m_j)\geq\e$.
Set $l_{j_0}=\max\{l_j:1\leq j\leq n\}$ and $f^n=f_{j_0}$.
Observe that for $1\leq j\leq n$, $f^n(m_{n+i})\geq\de$ for all
$i=1,\ldots,l_j$ and hence
$f^n(m_j)\geq\e$ for all $j=1,\ldots n$.
It is clear now that the \sq\ $\{f^n\}$ does not have weakly convergent
sub\sq. Therefore the case $[M]\sbs[N']\stm S$ is impossible
and it is easy to check that if $[M]\sbs S$ then $M$ satisfies
the conclusion of the lemma.

\th{1.2. Lemma} Consider $F$ relatively weakly compact subset of $c_0[\NN]$,
$\de>0$ and $0<\e<1$. Then \fe $N'\in[\NN]$ \txs
$M=(m_i)$ such that:

For every $f\in F$, $\nin$, $I\sbs\{1,\ldots,n\}$ with
$\min_{i\in I}f(m_i)\geq\de$ \txs $g\in F$ satisfying the
following two conditions

(i) $\min_{i\in I}g(m_i)>(1-\e)\,\de$

(ii) $\sum_{\{i\leq n:i\notin I\}}|g(m_i)|<\e\cdot\de$.

\proof Choose $a>0$ with $|f(m)|\leq a$ for $m\in\NN$, $f\in F$.
Next we choose a strictly increasing \sq\ $(k_n)$ of natural numbers
\st $2^{k_1}>a$ and if $\e_n=\frac{1}{2^{k_n}}$ then
$\sumnf\sum_{k=n}^\ty\e_k<\e\cdot\de$.

We divide the proof into two stages. In the first we will construct the
set $M$ and in the second we will show that it satisfies the
conclusion of the lemma.

\sprh{\bf The set $M=(m_i)$}

\vspace{4pt} The set $M=(m_i)$ is defined inductively so that the
following conditions are fulfilled.

If $I$ is a finite subset of \NN, and $j<\min I$, then \fe
$f\in F$ \st $\min_If(m_i)>\de$ \txs $g\in F$ with:

a) $\min_Ig(m_i)>\de$

b) $|g(m_j)|<\e_j$

c) $|g(m_i)-f(m_i)|\leq\e_j$ for $i=1,2,\ldots,j-1$.

To find such an $M$ we choose inductively a decreasing \sq\ of
infinite sets $N'\sps N_1\sps\cdots\sps N_i\sps\cdots$ and we
set $m_i=\min N_i$.

To choose $N_1$ we apply Lemma 1.1 to find $N_1$ subset of $N'$ \st
the conclusion of Lemma 1 holds for the given $\de$ and
$\e=\e_1$. This finishes the choice of $N_1$.

Suppose that $N'\sps N_1\sps\cdots\sps N_j$ have been chosen \st
if\\ 
$m_i=\min N_i$ then $m_1<m_2<\cdots m_j$ and if $1<i\leq j$,
$I$ is a finite subset of $N_i$ with $m_i<\min I$, 
$f\in F$ with $\min_{k\in I}f(k)>\de$ then \txs $g\in F$
satisfying (a), (b), (c).
To choose $N_{j+1}$ we consider the set $W$ of all closed
dyadic intervals of length $\e_{j+1}=\frac{1}{2^{k_{j+1}}}$
which are contained in the interval $[-2^{k_1},2^{k_1}]$.
We denote by $W^j$ the $j$-times product of $W$ and \fe
$\calb\in W^j$, $\calb=(B_1,\ldots,B_j)$, we set
\[ F_{\calb}=\{f\in F: f(m_i)\in B_i\},\]
which clearly is relatively weakly compact.

Applying repeatedly Lemma 1.1, we find a set $N_{j+1}$ infinite
subset of $N_j$ \st $m_j<\min N_{j+1}$ and the conclusion of
Lemma 1.1 holds for the set $N_{j+1}$ and \fe $F_{\calb}$,
$\calb\in W^j$, the given $\de$ and $\e=\e_{j+1}$.
This completes the inductive construction of the sets $(N_i)$ and
hence the set $M$ is defined.

It remains to show that $M$ satisfies the desired properties.

\sprh{\bf The set $M$ satisfies (I) and (II)}

\sprh Given $\nin$, $I$ subset of $\{1,2,\ldots,n\}$ and $f\in F$
\st\\
 $\min_{i\in I}f(m_i)>\de$, we shall define the desired
function $g$. For this, we inductively choose $g_0,g_1,\ldots, g_n$
elements of $F$ such that:

 $f=g_0$, 

if $k\in\nn$, $1\leq k<n$ and
$g_0,g_1,\ldots,g_k$ have been chosen satisfying
the property:

  \fe $1\leq l\leq k$ and $i=1,2,\ldots, l$, 

$|g_l(m_i)-g_{l-1}(m_i)|\leq\e_l$ and

 for $i\in\{l+1,\ldots,n\}\cap I$, we have that
$g_l(m_i)>\de$

and $g_l(m_l)>\de$ if $l\in I$ or $|g_l(m_l)|<\e_l$ otherwise.

To choose $g_{k+1}$ we distinguish two cases.

\noindent {\sl Case 1:} $k+1\in I$. Then we set $g_{k+1}=g_k$.

\noindent {\sl Case 2:} $k+1\notin I$. Then we choose $g_{k+1}$ \st
$|g_{k+1}(m_i)-g_k(m_i)|\leq\e_{k+1}$, $g_{k+1}(m_i)>\de$ \fe
$i\in\{k+2,\ldots,n\}\cap I$, $|g_{k+1}(m_{k+1})|<\e_{k+1}$.
The existence of such a $g_{k+1}$ follows from the properties
of the set $M$.

This completes the inductive definition of $g_0,\ldots, g_n$.
It is easy to see that the final function $g_n$ is the desired
$g$. The proof of the lemma is complete.

\th{1.3. Theorem} Every $(x_n)_{n\in{\bf N}}$ normalized weakly null \sq\
in a Banach space $X$ has a convexly unconditional sub\sq.

\proof Assume, by passing to a sub\sq\ if it is needed, that
 $(x_n)_{n\in{\bf N}}$ is Schauder basic with basic constant
$D\geq1$. We, inductively, apply Lemma 1.2 to choose a decreasing \sq\
 $(M_n)_{n\in{\bf N}}$ \st $M_n$ satisfies the conclusion of the
Lemma for $F=\{(x^*(x_n))_{n\in{\bf N}}$, $\|x^*\|\leq1\}$,
$\de=\frac{1}{n}$, $\e=\frac{1}{n^3}$.

We select a strictly increasing \sq\ $M= (m_n)_{n\in{\bf N}}$ \st\\
$m_n\in M_n$.

\noindent{\sl Claim:} The \sq\ $(x_n)_{n\in M}$ is
convexly unconditional. 

Indeed, given $x=\sum_{n\in M}a_nx_n$ an absolutely convex combination
with $\|x\|>\frac{1}{k}$ and $(\varepsilon_n)_{n\in M}$ a \sq\
of signs we choose $x^*\in B_{X^*}$ with $x^*(x)>\frac{1}{k}$.
There exists finite $J\sbs M$ \st
$x^*(\sum_{n\in J}a_nx_n)>\frac{1}{k}$.
We set
\[J_1=\left\{n\in J:|x^*(x_n)|>\frac{1}{2k}\right\}.\]
Then we have $x^*(\sum_{n\in J\stm J_1}a_nx_n)\leq\frac{1}{2k}$ and
hence
 $x^*(\sum_{n\in J_1}a_nx_n)>\frac{1}{2k}$.
By splitting the set
 $J_1$ into four sets, in the obvious way,  we find a subset
$I\sbs J_1$ such that:

$|\sum_{n\in I}\varepsilon_na_n|>\frac{1}{8k},$

$\{\varepsilon_na_n:n\in I\}$ are either all non-negative or all
negative and 

$\{x^*(x_n):n\in I\}$ are of the same sign.

\noindent We consider $x^*$ if the sign of $x^*(x_n)$ 
is positive and $-x^*$ if the
sign is negative and this we again denote by $x^*$.
Thus we also have $x^*(x_n)>\frac{1}{2k}$ for $n\in I$.
For every $r\in\NN$ we denote by $B(r)$ the unconditional constant
of the \ $\{x_{m_1},\ldots, x_{m_r}\}$. This means
that for $G\sbs\{1,\ldots,r\}$,
\[\left\|\sum_{i\in I}b_ix_{m_i}\right\|\leq B(r)\cdot
\left\|\sum_{i=1}^rb_ix_{m_i}\right\|.\]
(This happens because the norm $\|\cdot\|$ in the space of
dimension $r$ that is generated by $x_{m_1},\ldots, x_{m_r}$ is
equivalent to the maximum norm \wrt to this basis).

We split $I$ into two sets $I_1=I\cap\{m_1,\ldots,
m_{2k-1}\}$ and $I_2=I\stm I_1$. 

\bce{ We have
$|\sum_{n\in I_1}\vare_na_n|>\frac{1}{16k}$ or
$|\sum_{n\in I_2}\vare_na_n|>\frac{1}{16k}$.}

\noindent If the first condition holds, then

\beq\left\|\sum_{n\in M}\vare_na_nx_n\right\|&\!\!\geq\!\!&
\frac{1}{D}
\left\|\sum_{i=1}^{2k-1}\vare_{m_i}a_{m_i}x_{m_i}\right\|\geq
\frac{1}{D\cdot B(2k-1)}\left\|\sum_{n\in I_1}\vare_na_nx_n\right\|\\
&\!\!>\!\!&\frac{1}{D\cdot B(2k-1)}\left|\sum_{n\in I_1}
\vare_na_n\right|\cdot
\frac{1}{2k}>\frac{1}{D\cdot B(2k-1)}\cdot\frac{1}{32k}\rule{15pt}{0mm}(1)\eeq
In the  second case, \txs $y^*\in B_{X^*}$ \st

(i) $\min_{I_2}\{y^*(x_n)\}>\left(1-\frac{1}{(4k)^3}\right)\cdot
\frac{1}{2k}$ and

(ii) $\max\left\{|y^*(x_n)|:n\in\{m_{2k},\ldots,m_{l}\}\stm I_2\right\}<
\frac{1}{2k(4k)^3}$, where $m_l=\max I$.

\noindent Therefore
\beq \left|y^*\left(\sum_{i=2k}^l\vare_{m_i}a_{m_i}x_{m_i}\right)\right|&>&
\left|y^*\left(\sum_{n\in I_2}\vare_{n}a_{n}x_{n}\right)\right|
-\left|y^*\left(\sum_{n\in\{m_{2k},\ldots,m_l\}\stm I_2}
\vare_{n}a_{n}x_{n}\right)\right|\\
&>&\frac{1}{16k}\left(1-\frac{1}{(4k)^3}\right)\frac{1}{2k}-
\frac{1}{2k}\frac{1}{(4k)^3}\\
&=&
\frac{1}{32k^2}\left(1-\frac{1}{(4k)^3}-\frac{1}{(2k)^2}\right)>
\frac{1}{64k^2}\eeq

Finally

\ce{\left\|\sum_{n\in M}\vare_na_nx_n\right\|>\frac{1}{2D}
\left\|\sum_{i=2k}^l\vare_{m_i}a_{m_i}x_{m_i}\right\|>
\frac{1}{D\cdot 128k^2}}{(2)}

\noindent From (1) and (2) we get that
\[ C\left(\frac{1}{k}\right)\geq\min\left\{\frac{1}{D\cdot B(2k-1)}
\cdot\frac{1}{32k}\:,
\frac{1}{D\cdot 128k^2}\right\}.\]

\th{1.4. Corollary} Let  $(x_n)_{\nin}$ be a normalized weakly null
\sq. Then \fe $M\in[\nn]$ \txs $L\in[M]$, $L=(l_j)_{j\in\nn}$ \st
the following property is satisfied\\
For every $k>0$ \txs $(C(k)>0$ \st \fe $x=\sum_{j=1}^\ty a_jx_{l_j}$
and $\sum_{j=1}^\ty|a_j|<k$ then \fe  \sq\ $(\vare_j)_{j\in\nn}$
of signs we have
$C(k)<\lf\|\sum_{j=1}^\ty\vare_ja_jx_{l_j}\r\|\leq k$.

\proof Set $d=\sum_{j=1}^\ty|a_j|$ and $b_j=\f{a_j}{d}$.
Then $\sum_{j=1}^\ty b_jx_{l_j}$ is absolutely convex combination
and $\|\sum b_jx_{l_j}\|>\f{1}{d}>\f{1}{k}$. Hence
by Theorem 1.3 we get the left inequality with
$C(k)=C\lf(\f{1}{k}\r)$. The right is immediate from the
triangle inequality.

\sprh The following con\sq\ of Lemma 1.2 has been proved by
H. Rosenthal with the use of transfinite induction.

\th{1.5. Theorem} Let $K$ be a compact space and $(f_n)_{\nin}$
a \sq\ of continuous characteristic functions converging
pointwise to zero.  Then \txs $L\in[\nn]$, $L=(l_j)_{j\in\nn}$
\st  $(f_{l_j})_{j\in\nn}$ is an unconditional basic \sq.

\proof Define $F:K\lra c_0(\nn)$ by the rule
$F(x)=(f_n(x))_{\nin}$. Then for each $x\in K$,
$F(x)$ is a finite subset of $\nn$ and $F[K]$ is weakly compact.
By Lemma 1.2, \txs $L\in[\nn]$ \st \fe $x\in K$, $G\sbs
F(x)\cap L$ \txs $y\in K$ \st $F(y)\cap L=G$. It is
easy to check that the \sq\ $(f_n)_{n\in L}$ is an
unconditional \sq.

\sprt
\bce{{\large\bf 2. Summability methods}}

\sprh
{\bf Schreier Hierarchy, The RA Hierarchy}

\th{Notation} We denote by $S^+_{\ell^1}$ the positive part
of the unit sphere of $\ell^1(\NN)$. For $A=(a_n)_{n\in{\bf N}}$ in
$S^+_{\ell^1}$ and $F=(x_n)_{n\in{\bf N}}$ bounded \sq\ in
a Banach space $X$ we denote by $A\cdot F$ the usual matrices
product, that is:

\bce{$A\cdot F=\sumnf a_nx_n$.}

\th{2.1.1. Definition} For an $M$ infinite subset of \NN\ an
{\em $M$ summability method\/} is a block \sq\ $(A_n)_{\nin}$
with $A_n\in S^+_{\ell^1}$ and $M=\bigcupnf\spp A_n$
where $\spp A_n=\{\nin:a_n\neq 0\}$.

\th{2.1.2. Definition} Suppose that $(A_n)_{\nin}$ is an $M$ summability
method. A bounded \sq\ $F=(x_n)_{\nin}$ is said to be
{\em $M-(A_n)_{\nin}$ summable} if the \sq\ 
$(A_n\cdot F)_{\nin}$ is Cesaro summable. This means that the
\sq\ $z_n=\frac{\sum_{k=1}^nA_k\cdot F}{n}$ is norm convergent.

\th{2.1.3. Remark} To each $M\in[\NN]$, $M=(m_{n})_{n\in\NN}$,
we assign the $M$-summability method $A_n=\{e_{m_n}\}$. Then the
$M-(A_n)_{\nin}$ summability of $(x_n)_{\nin}$ is exactly the
usual Cesaro's summability of the sub\sq\ $(x_n)_{n\in M}$.

\sprh
{\bf Definition of Schreier Hierarchy}

\sprh Next we recall the definition of the generalized Schreier
families $(\calf_\xi)_{\xi<\omega_1}$. These are defined inductively
in the following manner. 

\th{Notation} For $F_1,F_2$ in $[\nn]^{<\om}$ we denote by 
$F_1<F_2$ the relation $\max F_1<\min F_2$.

\sprh
We set $\calf_0=\{\{n\}:\nin\}\cup\{\emptyset\}$.

Suppose that $\xi=\zeta+1$ and $\calf_\zeta$ has been defined. We set

$\calf_\xi=\left\{ F\in[\NN]^{<\omega}:F=\bigcup_{i=1}^nF_i,\;
F_i\in\calf_\zeta,\; n\leq F_1<\cdots<F_n\right\}$.

If $\xi$ is a limit ordinal and $\calf_\zeta$ has been defined for all
$\zeta<\xi$ then we fix a strictly increasing family of non-limit ordinals
$(\xi_n)_{\nin}$ with $\sup\xi_n=\xi$ and we define

$\calf_\xi=\left\{ F\in[\NN]^{<\omega}:n\leq\min F\;\:{\rm and}\;\:
F\in\calf_{\xi_n}\right\}.$

\th{Remark} The use of a \sq\ of non limit ordinals 
$(\xi_n)_{\nin}$ in the definition of $\cf_\xi$, $\xi$ limit,
is not important. We make this assumption in order to avoid
some more complexity of the Approximation Lemma  given below.

\sprh
{\bf Definition of the RA Hierarchy}

\sprh To each $M\in[\NN]$ and $\xi<\omega_1$ we will assign inductively
an $M$ summability method $(\xii^M_n)_{\nin}$ in the following manner.

(i) For $\xi=0,\;\: M=(m_n)_{\nin}$ we set $\xii_n^M=\{e_{m_n}\}$.

(ii)  If $\xi=\zeta+1$, $M\in\N$ and $(\zi^M_n)_{\nin}$ has been
defined then we, inductively, define
$(\xii^M_n)_{\nin}$ as it follows.
We  set $k_1=0$, $s_1=\min\spp\zi^M_1$, and
\[\xii^M_1=\frac{\zii^M_1+\cdots +\zii^M_{s_1}}{s_1}.\]
Suppose that for $j=1,2,\ldots,n-1$ $k_j,s_j$ have been defined and
\[\xii^M_j=\frac{\zii^M_{k_j+1}+\cdots +\zii^M_{k_j+s_j}}{s_j}.\]
Then  we set,
\[k_n=k_{n-1}+s_{n-1},\;\; s_n=\min\spp\zi^M_{k_n+1}\;\;
{\rm and}\]
\[\xii^M_n=\frac{\zii^M_{k_n+1}+\cdots +\zii^M_{k_n+s_n}}{s_n}.\]
This completes the definition for successor ordinals.

(iii) If $\xi$ is a limit ordinal and if we
suppose that \fe\ $\zeta<\xi$, $M\in\N$ the
\sq\ $(\zi^M_n)_{\nin}$ has been defined, then we define
$(\xii^M_n)_{\nin}$ as it follows:

\noindent We denote by
 $(\xi_n)_{\nin}$ the strictly increasing \sq\ of ordinals
with $\sup\xi_n=\xi$  that defines the family $\calf_\xi$.

For $M=(m_k)_{k\in{\bf N}}$ we inductively define $M_1=M$, 
$n_1=m_1$,

\noindent  $M_2=\left\{ m_k:m_k\not\in\spp[\xii_{n_1}]^{M_1}_1\right\},$
$n_2=\min M_2$,

\noindent
$M_3=\left\{ m_k:m_k\not\in\spp[\xii_{n_2}]^{M_2}_1\right\},$
and $n_3=\min M_3$, and so on.

We set 
\[\xii^M_1=[\xii_{\mbox{$n$}_1}]^{M_1}_1,\;
\xii^M_2=[\xii_{\mbox{$n$}_2}]^{M_2}_1,\ldots,\; 
\xii^M_k=[\xii_{\mbox{$n$}_k}]^{M_k}_1,\ldots.\]
Hence $(\xii^M_n)_{\nin}$ has been defined.
This completes the definition of RA Hierarchy.

\sprh\centerline{\bf Properties of the two Hierarchies}

The following properties can be established inductively.

{\bf P.1}: For $\xi<\omega_1$, $M\in\N$ $(\xii^M_n)_{\nin}$ is an 
M-summability method, i.e., $(\xii^M_n)_{\nin}$ is a block \sq\
of elements of $S^+_{\ell^1}$ and $M=\bigcupnf\spp\xii^M_n$.

{\bf P.2}: For every $\xi<\omega_1,$ $M\in\N$,
$\nin$ $\spp\xii^M_n\in\calf_\xi$.

{\bf P.3}: For every $\xi<\omega_1$ and every $N,M\in\N$ \st

$\spp\xii^M_i=\spp\xii^N_i$ for $i=1,2,\ldots,k$ we have

$\xii^M_i=\xii^N_i$ for $i=1,2,\ldots,k$.

{\bf P.4}: For every $M\in[\nn]$ and
$(n_k)_{k\in\nn}$ subset of $\nn$ if
$M'=\bigcup_{k=1}^{\infty}\spp\xii^M_{n_k}$ then
$\xii^{M'}_k=\xii^M_{n_k}$.

\th{Remark} Properties P.3 and P.4 are important for our proofs
and they indicate a strong stability of the methods
$(\xii^M_n)_{\nin}$.

\th{2.1.4. Definition} A family \calf\ of finite subsets of $\NN$
is said to be {\em adequate} if \calf\ is compact and \fe\ 
$F\in\calf$, if $G\sbs F$ then $G\in\calf$.

\th{2.1.5. Notation} If \calf\ is an adequate family and $L\in\N$ we set
$\calf[L]=\{F\in\calf:F\sbs L\}$.

Clearly $\calf[L]$ is an adequate subfamily of \calf.

\th{2.1.6. Notation} For an ordinal $\xi<\om_1$ and $M\in\N$, 
$M=(m_i)_{i\in{\bf N}}$ we define
\[\cf^M_\xi=\{G:G=(m_i)_{i\in F},\: F\in \calf_\xi\}.\]
It is easy to see that $\cf^M_\xi$ is an adequate family.

\th{2.1.7. Remark} It is proved readily by induction that
$\cf^M_\xi$ is a subfamily of $\calf_\xi[M]$; on the other
hand, it is not true that $\calf_\xi[M]$ is contained in $\cf^M_\xi$.
We will show that by going to a subset $N$ of $M$,
$\cf^M_\xi\N$ and $\calf_\xi\N$ are in a sense comparable.

\th{2.1.8. Lemma} For all ordinals $\z<\xi<\om_1$ \txs $n(\z,\xi)\in\NN$ \st
\fe $F\in\calf_1$, $n(\z,\xi)\leq\min F$, $F\in\calf_\xi$.

(b) The same holds for $\cf^M_{\zeta},\:\cf^M_\xi$.

\sprh The proof of this lemma is obtained easily by induction.

\th{2.1.9. Lemma} For $M\in\N$, $\xi$ a limit ordinal, $\e>0$ and $N\in[M]$
\txs $L\in\N$ satisfying the following property:

For every $P\in[L]$, $\nin$ \txs $G\in\cf^M_\xi$ \st
$\langle\xii^P_n,G\rangle>1-\e$.

\proof We proceed by induction. We will establish the following.

\sprh
\noindent{\bf Inductive hypothesis:} For every limit ordinal $\xi<\om_1$,
$N\in[M]$, $\e>0$ \txs $L\in\N$ such that:

For every ordinal $\z$, $\z\leq\xi$ \txs $l(\z,\xi)$ such that:

$\rule{1cm}{0mm}$For every $P\in[L]$, $\nin$ with $l(\z,\xi)\leq\min\spp
\zi^P_n$ \txs

$\rule{1cm}{0mm}$ $G\in\cf^M_\xi$ \st $\langle\zi^P_n,G\rangle>1-\e$.

\vspace{10pt}{\bf Note:} In the sequel for $M\in\N$, $\z<\xi$ we shall denote
by $k(\z,\xi)$ the natural number appearing in part (b) of the
previous lemma and satisfying the property:
If $F\in\cf^M_\z$, $k(\z,\xi)\leq\min F$ then $F\in\cf^M_\xi$.

\vspace{10pt}
We pass to prove the inductive hypothesis.

{\sl Case 1:} $\xi=0$. The proof follows immediately from the definitions.

\vspace{3pt}{\sl Case 2:} $\xi$ is a limit ordinal and \txs
an increasing \sq\ $(\z_k)_{k\in{\bf N}}$ of smaller limit ordinals
with $\sup\z_k=\xi$.
We also denote by $(\xi_n)_{\nin}$ the \sq\ of non limit ordinals, used
in the definition of the family $\calf_\xi$.

Let $N$ be a subset of $M$ and $\e>0$
We, inductively, choose a decreasing \sq\ $(L_k)_{k\in{\bf N}}$ of
subsets of $N$ \st each $L_k$ satisfies the inductive hypothesis for
the ordinal $\z_k$ and the number $\e'=\frac{\e}{4}$. Further,
for $\nin$ denote by $k_n$ the least natural number \st
$\xi_n<\z_{k_n}$. 

Inductively define a subset $L$ of $N$,
$L=(l_i)_{i\in{\bf N}}$ in the following manner.
Set  $M=(m_s)_{s\in{\bf N}}$, choose
$l_1\in L_1$ \st  $l_1=m_{s_0}$ and  $\frac{1}{2^{s_0}}<\frac{\e}{4}$.

In general choose $l_{i+1}\in L_{i+1}$ \st

$l_{i+1}>\max\{ l_i, k(\z_{k_{l_i}+2},\xi),\;
l(\xi^{-}_{l_i},\z_{k_{l_i}})\}$,

\noindent where $\xi^{-}_{l_i}$ denotes the predecessor of
$\xi_{l_i}$. 

The inductive definition of $L$ is complete and we prove the following:

\noindent{\sl Claim:} The set $L$ satisfies the inductive hypothesis
for the ordinal $\xi$ and the number $\e$.

Suppose first that $\z<\xi$. Then \txs $\z_k$ \st $\z<\z_k$. By the
inductive assumption \txs $l(\z,\z_k)\in\nn$ \st for $P\in[L_k]$,
$\nin$ with $l(\z,\z_k)<\min\spp\zi^P_n$ \txs $G\in\cf^M_{\z_k}$
with $\langle\zi^P_n,G\rangle>1-\e'$. We set
\[ l(\z,\xi)=\max\{l(\z,\z_k),\; k(\z_k,\xi)\}\]
and we show that it satisfies the inductive hypothesis for the pair
$(\z,\xi)$ and the number $\e$.

Indeed, if $P\in[L]$, $\nin$, with
$l(\z,\xi)\leq\min\spp\zi^P_n$ then \txs $G\in\cf^M_{\z_k}$ with
$ \langle\zi^P_n,G\rangle>1-\e'$.

Since $\cf^M_{\z_k}$ is adequate,
we assume that $l(\z,\xi)\leq G$. Then $k(\z_k,\xi)\leq G$
hence $G\in\cf^M_\xi$ and this completes the proof of this case.

We pass to the remaining case $\z=\xi$.

 We shall show that $l(\xi,\xi)=l_1$.

Indeed, suppose that $P\in[L]$, $\nin$. Then if $l_i=\min\spp
\xii^P_n$ then \txs $P_1\in[P]$ with
$\xii^P_n=[\xii_{\mbox{$l$}_i}]^{P_1}_1.$
Further $[\xii_{\mbox{$l$}_i}]^{P_1}_1=
\frac{1}{l_i}\sum_{j=1}^{l_i}
[\xii \mbox{$^{-}_{l_i}$}]^{P_1}_j.$ 

Set $F_j=\spp[\xii\mbox{$^{-}_{l_i}$}]^{P_1}_j$
and observe that:

(i)\hfill$F_1<F_2<\cdots <F_{l_i}$,\hfill\mbox{}

(ii)\hfill$l(\xi^-_{l_i},\z_{k_{l_i}})\leq l_{i+1}\leq F_2$\hfill\mbox{}

\noindent Hence by the inductive assumption \tx $G_2,G_3,\ldots,G_{l_i}$
in $\cf^M_{\z_{k_{l_i}}}$ \st $G_j\sbs F_j$ and
$\langle[\xii\mbox{$^-_{l_i}$}]^{P_1}_j,G_j\rangle>1-\e'$.

Notice also that $G_2<\cdots<G_{l_i}$.

We set $d_1=\left[\frac{l_i}{2}\right]$ and observe that if
$m_i=\min G_{d_1+1}$ then $i\geq d_1$ hence the set
$W_1=G_{d_1+1}\cup\cdots\cup G_{l_i}$ belongs to $\cf^M_{\z_{k_{l_i}}+1}$.
Repeat the same for $d_2=
\left[\frac{d_1}{2}\right]$ and for the set
$W_2=G_{d_2+1}\cup\cdots\cup G_{d_1}$ that also 
belongs to $\cf^M_{\z_{k_{l_i}}+1}$.

Therefore, following this precedure, we define $W_1, W_2,\ldots,
W_s$

 where  $s\in\nn$ \st $l_i=m_s$.

Since $l_1=m_{s_0}$,
$\frac{1}{2^{s_0}}<\e'=\frac{\e}{4}$ we get that
$\frac{1}{2^s}<\e'.$ 

Notice that
$m_s\leq W_s<W_{s-1}<\cdots<W_1$, hence the set
$G=\bigcup_{q=1}^sW_q$ belongs to  $\cf^M_{\z_{k_{l_i}}+2}$.

Further $\bigcup_{q=1}^sW_q=\bigcup_{j=j_0+1}^{P_i}G_j$ where
$j_0$ satisfies the following:

\bce{$j_0\leq l_i-\sum_{p=1}^s\frac{l_i}{2^p}=\frac{l_i}{2^s}<l_i\cdot\e'$.}
Finally
\[\lf\langle\xii^P_n,G\r\rangle =
\frac{1}{l_i}\sum_{j=j_0+1}^{l_i}\left\langle [\xii\mbox{
$^-_{l_i}$}]^{P_1}_j,
G_j\right\rangle
>(1-\e')-\frac{j_0}{l_i}>1-2\e'>1-\e.\]
Notice that  $k(\z_{k_{l_i}+2},\xi)\leq l_{i+1}\leq G$
and hence $G\in\cf^M_\xi$. 

This completes the proof for case 2.

\vspace{5pt}\noindent{\sl Case 3} $\xi=\z+\om$ with $\z$ a limit ordinal.

Let $N$ be a subset of $M$ and $\e>0$. By the inductive assumtion,
\txs $L\in\N$ satisfying the inductive hypothesis for the ordinal $\z$
and the number $\e'=\frac{\e}{4}$. Since every $L'\in[L]$
also satisfies the inductive hypothesis for the same $\z$ and $\e'$,
we may assume that $L$ satisfies the following.

If $M=(m_i)_{i\in\nn}$ and $L=(m_i)_{i\in D}$ then ordering $D$ in the
natural manner, we have that $L=(m_{i_n})_{\nin}$. Under this notation we
assume that 

\[\frac{1}{m_{i_1}}<\frac{\e}{4}\;\;{\rm  and}\;\;
\sum_{k=1}^\ty\frac{1}{2^{i_k}}<\frac{\e}{4}.\]

We first prove the following

\vspace{3pt}\noindent{\sl Claim 1:} For $n=1,2,\ldots$ \txs $k_n\geq n$
\st the following hold:

(i) If $P\in[L]$, $q\in\nn$ with $m_{i_n}\leq\spp[\zin]^P_q$,
\txs $G\in\cf^M_{\z+k_n}$ such that:

\bce{$\langle[\zin]^P_q,G\rangle >1-\left(\e'+\sum_{k=1}^n\frac{1}{2^{i_k}}
\right)>1-2\e'$.}

(ii) If $P\in[L]$, $q\in\nn$, $m_{i}=\min\spp[\zin]^P_q$ and

\bce{$[\zin]^P_q=\frac{1}{m_i}\sum_{j=1}^{m_i}
[\mbox{\boldmath$\zi+n-1$}
]^{P_1}_j$}

\noindent for some $P_1\in[P]$.

 If  $m_{i_n}\leq\spp[\mbox{\boldmath$\zi+n-1$}]^{P_1}_2$ then
\txs $G\in\cf^M_{j+k_n}$ \st

$m_{i_n}\leq\min G$ and
$\langle[\zin]^P_q,G\rangle >1-\left(2\e'+\sum_{k=1}^n\frac{1}{2^{i_k}}
\right)>1-3\e'$.

\vspace{3pt}\noindent{\sl Proof of Claim 1:}
We proceed by induction on $\nn$. The case $n=1$ is proved by similar
arguments as the general case.
Suppose that the claim has been proved for $k=1,2,\ldots, n-1$.
We prove first part (i).

Suppose that $P\in[L]$ \st 
$m_{i_n}\leq\spp[\zin]^P_1$. 

The case $[\zin]^P_q$ follows from this since there always 
exists $P_1\in[P]$ \st

\bce{$[\zin]^P_q=[\zin]^{P_1}_1$. }

We set
$m_i=\min\spp[\z+n]^P_1$ and write

$[\zin]^P_1=\frac{1}{m_i}\sum_{j=1}^{m_i}[\zinn]^P_j$.

 Also set
$F_j=\spp[\mbox{\boldmath$\zi+n-1$}]^P_j$
and notice that 

\mbox{}\hfill$F_1<F_2<\cdots<F_{m_i}$\hfill$(*).$

\noindent Since $m_{i_{n-1}}<F_j$, by the inductive assumption 
part (i), 
 \tx\\
 $G_j\in\cf^M_{\z+k_{n-1}}$, \st
$G_j\sbs F_j$ and
\[\langle[\mbox{\boldmath$\zi+n-1$}]^P_j,G\rangle >1-
\left(\e'+\sum_{k=1}^n\frac{1}{2^{i_k}}
\right).\] 
>From $(*)$ we also have,
$m_{i_n}\leq G_1<G_2<\cdots<G_{m_i}$. 

Set
$d_1=\lf[\f{m_i}{2}\r]$ and notice that 
 $m_{d_1+1}\leq\min G_{d_1}$ hence the set
$W_1=G_{d_1+1}\cup\cdots\cup G_{m_i}$ belongs to 
$\cf^M_{\z+k_{n-1}+1}$.
In the same manner we define $d_2=\lf[\f{d_1}{2}\r]$
and $W_2=G_{d_2+1}\cup\cdots\cup G_{d_1}$ that also belongs to 
$\cf^M_{\z+k_{n-1}+1}$.

Thus, successively define $W_1,W_2,\ldots,W_{i_n}$ with
$W_s\in\cf^M_{\z+k_{n-1}+1}$ for  $s=1,2,\ldots,i_n$ and

\bce{$m_{i_n}\leq W_{i_n}<W_{i_{n-1}}<\cdots<W_1$.}

\noindent Hence the set
$G=\bigcup_{s=1}^{i_n}W_s$ belongs to $\cf^M_{\z+k_{n-1}+2}$.

Further, $G=\bigcup_{s=1}^{i_n}W_s=\bigcup_{j=j_0+1}^{m_i}G_j$ and
$j_0$ is estimated by
\[j_0\leq m_i-\sum_{s=1}^{i_n}\f{m_i}{2^s}=
\f{m_i}{2^{i_n}}.\]
Finally we have,
\beq\langle[\zin]^P_1,G\rangle &=&
\f{1}{m_i}\sum_{j=j_0+1}^{m_i}
\langle[\zinn]^P_j,G_j\rangle\\
&>&1-\left(\e'+\sum_{k=1}^{n-1}\frac{1}{2^{i_k}}\right)-\f{j_0}{m_i}>
1-\left(\e'+\sum_{k=1}^{n}\frac{1}{2^{i_k}}\right)\eeq
and the proof of part (i) is complete.

The proof for part (ii)  is similar to the above one.

Indeed, for $P\in[L]$ we have as before

\bce{$[\zin]^P_1 =\f{1}{m_i}\sum_{j=1}^{m_i}
[\zinn]^P_j$}

\noindent  and if $F_j=\spp[\zinn]^P_j$
then by our assumption we have that $m_{i_n}\leq\min F_2$. 

Observe that
$[\zinn]^P_j$ for $j=2,\ldots,m_i$ satisfy part (i) of the inductive
assumption, hence \tx $G_j\sbs F_j$ \st
$G_j\in\cf^M_{\z+n_{k-1}}$ and
\[\lf\langle[\zinn]^P_j,G_j\r\rangle >
 1-\left(\e'+\sum_{k=1}^{n-1}\frac{1}{2^{i_k}}\right).\]
As in the previous part of the proof, we build $G$ in
$\cf^M_{n_{k-1}+2}$ \st
\[\lf\langle\f{1}{m_i}\sum_{j=2}^{m_i}
[\zinn]^P_j,G\r\rangle
>1-\left(\e'+\sum_{k=1}^{n}\frac{1}{2^{i_k}}\right).\]
Further, since $\f{1}{m_i}\leq\f{1}{m_1}<\e'$, we get that
\[\lf\langle [\zin]^P_1,G\r\rangle >
1-\left(2\e'+\sum_{k=1}^{n}\frac{1}{2^{i_k}}\right).\]
The proof of Claim 1 is complete.

To finish the proof of the lemma for case 3, it remains to define
an appropriate $L'$ subset of $L$.
The set $L'=(l_j)_{j\in\nn}$ is defined by induction.
We denote by $(\xi_k)_{k\in\nn}$ the increasing \sq\ of smaller than
$\xi$ ordinals \st $\lim\xi_k=\xi$. Notice that, since
$\xi=\z+\om$ \txs $k_0$ \st
for all $k\geq k_0$ $\xi_k=\z+n_k$.
Choose $l_1\in L$ with $l_1>k_0$. Suppose that
$l_1<\cdots<l_j$ have been defined. To define $l_{j+1}$ we follow
the next procedure. Since $l_j>l_1>k_0$, \txs $n_{k(j)}$ \st
$\xi_{l_j}=\z+n_{k(j)}$. From Claim 1 \txs $k(n_{k(j)})$
satisfying part (i) and (ii) of claim 1.
Further \txs $k'_j$ \st $\z+n_{k(j)}<\xi_{k'_j}$.
Choose $l_{j+1}\in L$ such that:
$l_{j+1}>\max\{ l_j,m_{i_j},k'_j\}$.
This completes  the inductive of $L'$.

\noindent{\sl Claim 2:} The set $L'$ satisfies the inductive hypothesis
for the ordinal $\xi$ and the number $\e$.

Indeed, given $\eta\leq\xi$

\noindent{\sl Step 1:} $\eta\leq\z$. Then using the inductive hypothesis
for $\z$ and $\eta$ we can easily establish the inductive
hypothesis for $\xi$ and $\eta$.

\noindent{\sl Step 2:} $\z<\eta<\xi$. Then $\eta=\z+n$ and we set
$l(\eta,\xi)=l_{j_\eta}$ and by definition
 $l_{j_\eta}>\max\{m_{i_n},k_\eta\}$
with $k_\eta\in\nn$ \st  $\xi_{k_\eta}>\z+k(n)$ and
$k(n)$ is the number corresponding to $n$ in Claim 1. 

If $P\in[L']$ and $n_0\in\nn$ and $l_{j_n}\leq\min\spp\etai^P_{n_0}$
then by Claim 1 \txs $G\in\cf^M_{\z+k_n}$ \st
$\langle\etai^P_{n_0},G\rangle >1-2\e'$.

Further, $G\in\cf^M_{\xi_{k_n}}$ and $\min G>k_n$ which
implies that $G\in\cf^M_\xi$ and the proof of Step 2 is complete.

\noindent{\sl Step 3:} $\eta=\xi$. In this case we shall show that
$l(\xi,\xi)=l_1$. 

Indeed, if $P\in[L']$ and $n_0\in\nn$ then,
by definition, for $l_j=\min\spp\xii^P_{n_0}$ choose
$P_1\in[P]$ such that:
$\xii^P_{n_0}=[\xii_{l_j}]^{P_1}_1$.

 Further,
\[   [\mbox{$\xii_{l_j}$}]^{P_1}_1
=\f{1}{l_j}\sum_{i=1}^{l_j}[\mbox{$\xii^-_{l_j}$}]^{P_1}_i\]
where $\xi^-_{l_j}$ denotes the predecessor of $\xi_{l_j}$.

Also, $\xi_{l_j}=\z+n_{k(j)}$ and $l_{j+1}>m_{i_j}$.
Hence $[\mbox{$\xii_{l_j}$}]^{P_1}_1$ 
satisfies part (ii) of Claim 1 and therefore \txs
$G\in\cf^M_{\xi_{k'_j}}$ with $\min G\geq l_{j+1}$ and
\[\lf\langle[\mbox{$\xii_{l_j}$}]^{P_1}_1,G\r\rangle 
>1-3\e'>1-\e.\]
Since $l_{j+1}\geq k'_j$, we get that $G\in\cf^M_\xi$ and
the proof of step 3, as well as the proof of Case 3 and the proof of
the lemma are complete.

\th{2.1.10. Proposition} Let
 $\xi<\om_1$. Then  \txs $\de_\xi>0$ satisfying the following
property\\
 For every $M\in[\nn]$ \txs $L\in[M]$ 
 \st \fe $P\in[L]$, $\nin$ \txs $G\in\cf^M_\xi$ with
$\langle\xii^P_n,G\rangle >\de_\xi$.

\proof If $\xi$ is a limit ordinal then we set $\de_\xi=\f{1}{2}$ and the
previous lemma proves the desired result.

If $\xi$ is a successor ordinal 
then  $\xi=\z+n$ with $\z$ a limit ordinal.
Then we inductively
prove that  $\de_\xi\geq\f{1}{2^{n+1}}$.

Indeed, if $\xi=\z+1$ choose $L\in[M]$ \st the conslusion of the
previous lemma is satisfied for the ordinal $\z$
and the number $\e=\f{1}{2}$. We show that
$L$ satisfies the conclusion of the proposition foor $\xi=\z+1$ and
$\de_\xi=\f{1}{4}.$ 
This is so, since for $P\in[L]$, $\nin$
\txs $P_1\in[P]$ \st $\xii^P_n=\xii^P_1$ and further, if
$l=\min\spp\xii^{P_1}_1$ then
$\xii^{P_1}_1=\f{1}{l}\sum_{j=1}^l\zi^{P_1}_j$.
Choose $G_j\sbs\spp\zi^{P_1}_j$ \st $G_j\in\cf^M_\z$ and
$\langle\zi^{P_1}_j,G_j\rangle >\f{1}{2}$ and it is easy
to see that if $d=\lf[\f{l}{2}\r]$ then the set
$G=G_{d+1}\cup\cdots\cup G_l$ is in $\cf^M_{\z+1}$ and that
$\langle\xii^{P_1}_1,G\rangle >\f{1}{4}$.
This completes the proof for $\xi=\z+1$. The general case is proved
by a similar argument. The proof is complete.

\th{2.1.12. Remark} A con\sq\  of the above Lemma is that
Schreier hierarchy is in a sense universal.

Indeed, consider $f:\nn\lra\nn$ any strictly increasing function and
define 
\[\calf^f_1=\{ F\sbs\nn:\min F=n,\:\# F\leq f(n)\}.\]
$\calf^f_1$ is an adequate family and the regular $\calf_1$ is
$\calf^f_1$ for $f={\rm id}_{\nn}$.

By iteration we produce $(\calf^f_\xi)_{\xi<\om_1}$ and the repeated
averages hierarchy $[\xii^f]^M_n$ for $M\in\N$, $\nin$.

Next define the set $M=(m_i)_{i\in\nn}$ by the rule
$m_i=f(i)$. Then observe that if

\bce{$\calf^{f,M}_\xi=\{(m_i)_{i\in F},\: F\in\calf^f_\xi\}$,}

\noindent we have that

\bce{$\calf^{f,M}_\xi=\calf_\xi[M]$.}

\noindent Therefore from the above lemma we get that  \txs
$L\in[M]$ \st \fe $P\in[L]$, $\nin$ \txs
$G\in\cf^M_\xi$ with $\langle [\xii^{f,M}]^P_n,G\rangle >\de_\xi$.

This shows that $\calf^{f,M}_\xi$, $\cf^M_\xi$ are comparable
on the set $L$.

\sprh{\bf 2.1.13. Lemma} {\sl (Approximation Lemma)} Let $\xi<\omega_1$,
$M\in\N$, $\epsilon>0$. We set 
 $W={\rm co}\,(\{\xi^N_n:\nin, N\in[M]\})$.\\ 
Then for every ordinal $\zeta$ \st
$\xi\leq\zeta<\omega_1$, $L\in[M]$ \txs $L_\zeta\in[L]$ satisfying
the following property:

For every  $L'\in[L_\zeta]$, $\nin$ we have that
\[d_{\ell^1}(\zeta^{L'}_n,W)<\epsilon.\]

\proof Fix $\xi<\om_1$ $M\in\N$. We shall prove it by induction
for $\z$ greater than $\xi$, every $L\in[M]$ and $\e>0$.

\noindent (i) $\z=\eta+1$. Indeed, if $M\in\N$, $\e>0$ then \tx
$L_\eta$ satisfying the conclusion for the ordinal $\eta$.
Set $L_{\z}=L_\eta$ and it obvious that \fe $L'\in[L_\z]$,
$\nin$ we get the desired property
$d_{\ell^1}(\zeta^{L'}_n,W)<\epsilon$.

\noindent (ii) $\z$ is a limit ordinal.
Then fix the strictly increasing \sq\ $(\z_n)_{\nin}$ of successor
ordinals \st $\sup\z_n=\z$ and $(\z_n)_{\nin}$ defines the
family $\calf_\z$. 
Since each $\z_n$ is a successor ordinal it has the form
$\z_n=\xi_n+1$. 

 Choose $L_0\in[M]$
with $\min L_0=m_1$ and $\frac{1}{m_1}<\frac{\e}{4}$.
We inductively choose $L_0\sps L_1\sps\cdots\sps L_k\sps\cdots$
\st if $n_k=\min L_k$ then $(n_k)_{k\in{\bf N}}$ is strictly increasing
and $L_k=L_{n_{k-1}}$ for $M=L_{k-1},$ $\z=\xi_{n_{k-1}}$, $\frac{\e}{2}$.

\noindent{\sl Claim:} The set $N=(n_k)_{k\in{\bf N}}$ is the desired
$L_\xi$.

Indeed, let $L'\in[L_\xi]$, $\nin$. Then, by definition,
$\zi^{L'}_k=[\mbox{$\zi_{n_k}$}]^{L'_k}_1$, where 
$n_k=\min\spp\zi^{L'}_n$, $L'_k=\{ m\in L':
n_k\leq m\}$. It is clear that $L'_k\stm\{n_k\}\sbs L_k$.
Since $\z_{n_k}=\xi_{n_k}+1$, again by definition,
$[\zi_{n_k}]^{L'_k}_{1}$ is an average of $n_k$ many successive elements of
$\lf( [\mbox{$\xii_{n_k}$}]^{L'_k}_{n}\r)_{\nin}$. 
Since all of them except
the first one are $\frac{\e}{2}$ approximated by convex combination
of $W$ and $\frac{1}{n_k}<\frac{\e}{4}$ we get that
$d_{\ell^1}([\mbox{$\zi_{n_k}$}]^{L'_k}_1,W)<\epsilon$ 
and hence $d_{\ell^1}(\zi^{L'}_k,W)<\epsilon$.
This completes the proof of the lemma.

\sprt{ \bf 2.2 Strong Cantor-Bendixson index}

\th{2.2.1 Definition} Let \calf\ be an adequate family. 
For $L\in[\nn]$ we define the
strong Cantor-Bendixson derivative of $\calf[L]$ by the rule:

$\calf[L]^{(1)}=\{A\in\calf[L]:\forall\: N\in[L]: \: A$ is a limit
point of $\calf[A\cup N]\}$.

\th{2.2.2 Remark} It is clear that $\calf[L]^{(1)}$ is closed and nowhere
dense subset of $\calf[L]$ which is also adequate.

\sprh If $\xi=\z+1$ then we, inductively, define
$\calf^{(\xi)}[L]=\left(\calf^{(\z)}[L]\right)^{(1)}$ 
if $\xi$ is a limit ordinal then we set
\[\calf^{(\xi)}[L]=\bigcap_{\z<\xi}\calf^{(\z)}[L].\]

We define the {\em S.C.B. index of\/} $\calf[L]$ as the smallest
ordinal $\xi_0$ \st $\calf^{(\xi_0)}[L]=\emptyset$.

We denote this index by $s(\calf[L])$.

\th{2.2.3 Proposition} If $\xi$ is a limit ordinal and $L\in\N$ \st\\
$ s(\calf[L])>\xi$ then  \fe $N\in[L]$ we have that
$s(\calf[N])>\xi$.

\proof We will show that \fe ordinal $\z$ satisfying $\z<\xi$
we have that $s(\calf[N])>\z$.

Indeed, if $\z<\xi$ then $\calf^{\z+1}[L]\neq\emptyset$ hence
\txs $A\in\calf[L]$ with $A$ a limit point of
$(\calf^\z[L])[A\cup N]$. 
Since $\calf^\z[L]$ is adequate, we get that \txs at  least one
nonempty subset of $N$ that belongs to $\calf^\z[L]$. It is also easy
to see that $\calf^\z[L]\cap [N]^{<\om}\sbs\calf^\z[N]$ hence
$\calf^\z[N]\neq\emptyset$ and $s(\calf[N])>\z$.
Since $\xi$ is a limit ordinal and $s(\cf[N])>\z$ for all ordinals
$\z<\xi$ we get that $s(\cf[N])>\xi$.

This completes the proof.

\th{2.2.4 Proposition} If \calf\ is  an adequate family, $\xi$ is a limit
ordinal, $L\in\N$, \st $L$ is almost contained in $N$, then
$s(\calf[N])>\xi$ implies that $s(\calf[L])>\xi$.

\proof Similar to the previous one.

\th{2.2.5 Notation} In the sequel we will denote by
$\calf[L]^\xi$ the $\xi$-derivative of $\calf[L]$, while for
$N\in[L]$ we denote by
$\calf^\xi[N]$ the restriction of
$\calf^\xi[L]$ on the set $N$.

\th{2.2.6. Theorem} Let \calf\ be an adequate family. If $L\in\N$ \st
 $s(\calf[L])>\om^\xi$ then \txs $M\in[L]$, $M=(m_i)_{i\in{\bf N}}$
satisfying the property:
$\cf^M_\xi$ is a subfamily of $\cf[M]$.

\proof We prove it by induction. For $\xi=0$ the result is obvious.
Suppose that we have proved it for all $\z<\xi$. To prove it for
$\xi$ we will use a method created by Kiriakouli -- Negrepontis [M-N].
This method consists in a double induction. We start with the
next definition.

\th{Definition} An $n$-tuple $(\xi_1,\ldots,\xi_n)$ has property
(A) if \fe adequate gamily
$\calf$ with 
\[s(\calf)>\om^{\xi_n}+\om+\cdots+\om^{\xi_1}+\om\]
and \fe $L\in\N$ \txs $N\in[L]$ \st $N=(n_i)_{i\in{\bf N}}$,
and \fe 
\[ F_1\in\calf_{\xi_1},\ldots,F_n\in\calf_{\xi_n}\;\;{\rm with}\;
F_1<F_2<\cdots<F_n\] 
the set
$\{n_i:i\in\bigcup_{k=1}^nF_k\}$ belongs to $\calf[N]$.

\th{2.2.7 Lemma} Suppose that $(\xi_1,\xi_2,\ldots,\xi_n)$
has property (A) and $\z<\om_1$. Then $(\z,\xi_1,\xi_2,\ldots,\xi_n)$
has also property (A).

\proof We prove, by induction that \fe $\z<\om_1$, $l\in\nn$ the tuple
\[ (\underbrace{\z,\ldots,\z}_{l{\rm-times}},\xi_1,\ldots,\xi_n)
\;{\rm has\; property\;} A.\]
{\sl Case 1.} $\z=0$. Given $(\eta_1,\ldots,\eta_k)$ with the property
(A) we show that $(0,\eta_1,\ldots,\eta_k)$ has also property (A). 

Indeed,
start with  \calf\ adequate \st 
\[s(\calf)>\om^{n_k}+\om+\cdots+\om^{n_1}+\om+1+\om.\]
Since $s(\calf)$ is greater than a limit ordinal we get that
for every $L\in[\nn]$
\[ s(\calf[L])>\om^{n_k}+\om+\cdots+\om^{n_1}+\om+\om.\]
Therefore if we denote by $\z=\om^{n_k}+\om+\cdots+\om^{n_1}+\om$
we get that $\calf[L]^\z$ is an infinite set and since $\calf[L]^\z$
is an adequate family, \txs $M=(m_i)_{i\in{\bf N}}$ \st
$\{\{m_i\}:i\in\nn\}$ is a subfamily of
$\calf[L]^\z$.

Observe that the set

\centerline{$\calg_{m_1}=\{ F\in\calf[M\stm\{m_1\}]:
\{m_1\}\cup F\in\calf[M]\}$}

\noindent is an adequate family and $s(\calg_{m_1})>\z$. 
This is so since
$s(\calf[M])>\z$ and $\{m_1\}\in\calf[M]^\z$. 
>From the inductive
assumption \txs $M_1\in[M]$ \st \fe
\[ F_1\in\calf_{n_1},\ldots,
F_k\in\calf_{n_k}\;\;{\rm with}\;\;
F_1<F_2<\cdots<F_k\]
implies that
\[D=\{m^1_i:i\in\bigcup_{j=1}^k\calf_{n_j}\}\in\calg_{m_1}[M_1],\]
where $M_1=(m^1_i)_{i\in{\bf N}}$. 

Then clearly for any such $D$ the set $\{m_1\}\cup D$ belongs to
$\calf[M]$.

Set $n_1=m_1$, $n_2=m^1_1$ and repeat the same procedure by
defining $\calg_{n_2}$ and finding
$M_2\in[M_1]$ \st $M_2=(m^2_i)$, $n_2<m^2_1$ and if 
\[ F_1\in\calf_{n_1},
\ldots, F_k\in\calf_{n_k}\;\;{\rm satisfying}\;\;
F_1<F_2<\cdots<F_k\]
then the set
$\{n_2\}\cup\{m^2_i:i\in\bigcup_{j=1}^kF_j\}$ belongs to $\calf[M_2]$.
Following the same procedure, we, inductively,
 choose $n_l,\: M_l$ with 
\[ M=M_0\sps M_1\sps M_2\sps\cdots\sps M_l,\;\;{\rm and}\;\;
n_l\in M_{l-1}\]
satisfying the above properties. It follows now immediately that
the set $N=(n_l)_{l\in{\bf N}}$ satisfies the required properties
and hence $(0,\eta_1,\eta_2,\ldots,\eta_k)$ has property (A).

\vspace{4pt}
{\sl Case 2.} $\xi=\z+1$. 
Then by the inductive hypothesis, \fe $k$-tuple
$(\eta_1,\eta_2,\ldots,\eta_k)$ with the property (A) and every
$l\in\NN$ the $l+k$-tuple\\
$(\underbrace{\z,\ldots,\z}_{l{\rm-times}},$
$\eta_1,\ldots,\eta_k)$
has the property (A).

For every  $L\in\N$ \st
\[ s(\calf[L])>\om^{n_k}+\om+\cdots+\om^{n_1}+\om+\om^\xi+\om\]
we have that \fe $l\in\nn$
\[ s(\calf[L])>\om^{n_k}+\om+\cdots+\om^{n_1}+\om+\underbrace{\om^\z+\om+
\cdots+\om^\z+\om}_{l{\rm-times}}.\] 
Hence we can find $L\sps L_1\sps \cdots\sps L_l\sps\cdots$ with
$L_l$ satisfying the property:

\vspace{-6pt}
\begin{quote}
if $F_1\in\calf_\z,\ldots,F_l\in\calf_\z,F_{l+1}\in\calf_{\eta_1},
\ldots,F_{l+k}\in\calf_{\eta_k}$\\
 and $F_1<\cdots<F_{l+k}$ then
$\left\{m^l_i:i\in\bigcup_{j=1}^{l+k}F_j\right\}\in\calf[L_l]$.
\end{quote}

\vspace{-6pt}
\noindent Then if $N=(n_l)_{l\in{\bf N}}$ with $n_l\in L_l$, it
is easy to see that $N$ satisfies the required properties hence
$(\xi,\eta_1,\ldots,\eta_k)$ has property (A).

\vspace{4pt}
\noindent{\sl Case 3.} $\xi$ is a limit ordinal.
The proof is similar to the previous case.

\vspace{6pt}{\bf Proof of the theorem:} 
 \noindent{\sl Case 1.} $\xi=\z+1$.

 Since $s(\calf)>\om^\xi>l\cdot
(\om^\z+\om)$ and $(\underbrace{\z,\z,\ldots,\z}_{l{\rm-times}})$
has property (A) then for any $L\in\N$ choose
\[ L=L_0\sps
L_1\sps L_2\sps\cdots\sps L_l\sps\cdots\] \st $L_l$
 witnesses
property (A) for the set $\calf[L_{l-1}]$ and the $l$-tuple 
$(\z,\ldots,\z)$.\\
It is easy to show
that for any $N\in[L]$ \st $N=(n_l)_{l\in{\bf N}}$ and
$n_l\in L_l$ $\calf[N]$ satisfies the inductive assumption
for the ordinal $\xi$.

\vspace{4pt}
\noindent{\sl Case 2.} $\xi$ is a limit ordinal.
Let $(\xi_n)_{\nin}$ be the strictly increasing \sq\ with
$\sup\xi_n=\xi$ that defines the family $\cf_\xi$.
 For every $L\in\N$ we choose
$L=L_0\sps L_1\sps\cdots\sps L_n\sps\cdots$ \st
$L_n$ witnesses property (A) for the $(\xi_n)$ and the set
$\calf[L_{n-1}]$. If we set $N=(k_n)_{\nin}$ \st $k_n\in L_n$
then we easily check that $\calf[N]$ satisfies the inductive
assumption for the ordinal $\xi$.
The proof of the theorem is complete.

\vspace{10pt}
{ \bf 2.3 Large families}

\vspace{7pt}{\bf 2.3.1 Definition:} Let \calf\ be an adequate family, 
$M\in\N$, $\xi<\om_1$ and $\e>0$. We say that \calf\ is
{\em $(M,\xi,\e)$ large} if \fe $N\in[M]$ and every $\nin$ we have that
\[\sup_{F\in{\cal F}}\langle\xii^N_n,F\rangle>\e.\]

\th{2.3.2 Proposition} If \calf\ is an adequate family
which is $(M,\xi,\e)$ large. Then \fe $L\in[M]$ \txs $N\in[L]$
\st $s(\calf[N])>\om^\xi$.

\sprh
This proposition is one of the basic ingredients for the proof of the main
Theorems of this section. The result of this in connection with
Theorem 2.2.6 shows that every $(M,\xi,\e)$ large family \calf\ contains
the family $\calf^{N'}_\xi$. Hence the summability methods 
$\{(\xii^N_n)_{\nin},N\in[M]\}$ are sufficiently many to describe the
Schreier family $\calf^{N'}_\xi$. The proof of the proposition depends
strongly on infinite Ramsey theorem and the stability properties
P.3 -- P.4 of the RA hierarchy.

\sprh{\bf Proof of the proposition:}
We proceed by induction. The inductive hypothesis is the statement
of the proposition.

\vspace{3pt}
\noindent{\sl Case 1.} $\xi=0$. This is the easiest case since the result
immediately follows from the definitions.

\vspace{4pt}
\noindent{\sl Case 2.} $\xi$ is a limit ordinal. In this case
we prove first the following.

{\sl Claim:} For every ordinal $\z$  with $\z<\xi$ 
and every  $L\in[M]$ \txs $N\in[L]$ \st
$s(\calf[N])>\om^\z$.

Indeed, given $L\in[M]$ we define a partition of $[L]$ into
$A_1, A_2$ by the rule:
\[ A_1=\left\{N\in[L]:\sup_{F\in{\cal F}}
\langle\zi^N_1,F\rangle\leq\frac{\e}{2}\right\},\] $A_2=[L]\stm A_1.$
Notice that if $N=(m_i)_{i\in{\bf N}}$ and $N'=(m'_i)_{i\in{\bf N}}$
are \st  $m_i=m'_i$ for all $i\leq k=\max\spp\zi^N_1$ then by
P.3 we get that $\zi^N_1=\zi^{N'}_1$ hence $A_1$ is an open set.
Therefore from infinite Ramsey theorem we get that \txs
$L_1\in[L]$ \st either $[L_1]\sbs A_1$ or $[L_1]\sbs A_2$.

Assume that $[L_1]\sbs A_1$. Then by P.4 we have that \fe
$N\in[L_1]$ and every $\nin$ we have that  $\sup_{F\in{\cal F}}
\langle\zi^N_n,F\rangle\leq\frac{\e}{2}$.\\
This is so since any such $\zi^N_n$ is equal to
$\zi^{N'}_1$ for some $N'\in[N]$.\\ 
But then, from Lemma 2.1.13
\txs $L_\xi\in[L_1]$ \st \fe $\nin$
$d_{\ell^1}(\xii^{L_\xi}_n,W)<\frac{\e}{2}$ where 
$W={\rm co}(\{\zi^N_n:\nin,\: N\in[L_1]\}$).
Hence $\sup_{F\in{\cal F}}\langle\xii^{L_\xi}_n,F\rangle\leq\frac{\e}{2}
+\frac{\e}{2}=\e$, a contradiction, therefore $[L_1]\sbs A_2$.

This means that the family $\calf[L_1]$ is 
$\left(L_1,\z,\frac{\e}{2}\right)$
large and by the inductive assumption we get that \txs $N\in[L_1]$
\st $s(\calf[N])>\om^\z$.
Next choose a strictly increasing \sq\ of ordinals $(\xi_n)_{\nin}$
with $\sup\xi_n=\xi$. 

Inductively we choose
$L=L_0\sps L_1\sps L_2\sps\cdots\sps L_n\sps\cdots$ \st
$s(\calf[L_n])>\om^{\xi_n}$.

It is easy to see that \fe $N$, $N\in[L]$ 
\st $N$ almost contained in $L_n$
we have that $s(\calf[N])>\om^{\xi_n}$ and therefore 
$s(\cf[N])>\om^\xi$. The proof for case 2 is complete.

\vspace{3pt}
\noindent{\sl Case 3.} $\xi=\z+1$. We start with the following Lemma,
the proof of which uses again infinite Ramsey Theorem.

\th{2.3.3. Lemma} Let $\xi=\z+1$, $M\in\N$, $\e>0$ and \calf\ 
an adequate family that is $(M,\xi,\e)$ large. Then \fe $L\in[M]$
and every $\nin$ \txs $L_n\in[L]$ \st \fe $N\in[L_n]$ and $k\in\NN$
\[\sup_{F\in{\cal F}}\min\{\langle\zi^N_k,F\rangle:1\leq k\leq n\}>
\frac{\e}{2}.\]

\proof Consider $L\in[M]$ and $\nin$, and define a partition of $[L]$
into $A_1,A_2$ by the rule
\[ A_1=\{N\in[L]:\exists\: F\in\calf,\{\langle\zi^N_k,F\rangle>
\frac{\e}{2}\;\;{\rm for}\; k=1,\ldots,n\}\]
 and 
$A_2=[L]\stm A_1$.

As in the previous lemma, $A_1$ is an open set hence by infinite Ramsey
Theorem \txs $L_n\in[L]$ \st $[L_n]\sbs A_1$ or $[L_n]\sbs A_2$.

We will show that the second case is not possible and this will
prove the lemma.

Indeed, assuming that $[L_n]\sbs A_2$ we get that \fe 
$N\in[L_n]$ and every $k_1<k_2<\cdots<k_n$ and every $F\in\calf$

\ce{\min\lf\{\lf\langle\zi^N_{k_i},F\r\rangle,i=1,2,\ldots n\r\}\leq
\frac{\e}{2}.}{(1)}

\noindent This follows from the fact that \txs $N'\sbs N$ \st $\zi^{N'}_i=
\zi^N_{k_i}$ for all $i=1,\ldots,n$.\\
Each $\xii^{L_n}_m$ is an average of successive elements of
$(\zi^{L_n}_k)_{k\in{\bf N}}$, that is
$\xii^{L_n}_m=\frac{1}{s_m}\lf(\zi^{L_n}_{k_m}+\cdots+
\zi^{L_n}_{k_m+s_m}\r)$ and $(s_m)_{m\in{\bf N}}$ is stictly
increasing.
Choose $F\in\calf$ \st
$\langle\xii^{L_n}_m,F\rangle >\e$.

Then for large $s_m$ we get that
$\#\left\{i:\langle\zi^{L_n}_{k_m+1},F\rangle >\frac{\e}{2}\right\}>n$.\\
But this contradicts (1) and the proof is complete.

\th{2.3.4 Lemma} Assume that $\xi, M,\e,\calf$ are as in the previous
lemma. Then \fe $L\in[M]$ \txs $N\in[L]$ \st \fe $N'\in[N]$ with
$\min N'\geq n$ we have that
\[\sup_{F\in{\cal F}}\min\{\langle\zi^{N'}_k,F\rangle :k=1,\ldots,n\}
>\frac{\e}{2}.\]

\proof We apply, inductively, the previous lemma and we  choose\\
$ L\sps L_1\sps\cdots\sps L_n\sps\cdots$ \st the set
$L_n$ satisfy
the conclusion of the previous lemma for the number $n$. 
Then any set $N=(m_n)_{\nin}$ with the
property $m_n\in L_n$ has the desired property.

\th{2.3.5 Lemma} Let $\z<\om_1$, $M\in\N$, $\e>0$ and \calf\ be an
adequate family. Suppose that for some $\nin$ we have that \fe
$L\in[M]$
\[\sup_{F\in{\cal F}}\min\{\langle\zi^{L}_k,F\rangle :k=1,\ldots,n\}
>\e.\]
Suppose that  $\z$ satisfies the inductive assumption.
Then \fe
$L\in[M]$ \txs $N\in[L]$ \st
\[ s(\calf[N])>n\cdot\om^\z.\]

\proof We proceed by induction on $\NN$.

\vspace{3pt}
\noindent{\sl Case 1.} $k=1$. As we have shown in previous proofs, the
fact that for $N\in[M]$ $\sup_{F\in{\cal F}}\langle\zi^N_1,F\rangle>\e$
implies that \calf\ is $(M,\z,\e)$ large hence by the inductive
assumption every $L\in[M]$ contains infinite subset $N$ \st\\
$ s(\calf[N])>\om^\z.$

\vspace{3pt}
\noindent{\sl Case 2.} $k=n$. Assume that the Lemma has been proved
for all $k=1,2,\ldots, n-1$.

Given $L\in[M]$, since \fe $N\in[L]$ the vector 
$\zi^N_1$ has finite support
and rational coefficients we get that the set
$\{\zi^N_1:N\in[L]\}$ is countable and we order it as
$(\zi_n)_{\nin}$.
Consider $\zi_1$ and fix $L_1\in[L]$ with 
\[\max\spp\zi_1<\min L_1.\]

Let $\{F_i\}_{i=1}^d$ be an enumeration of all nonempty subsets of
$\spp\zi_1$. We define a partition of $[L_1]$ into a family
$(A_i)_{i=1}^d$ defined by the rule

$A_i=\{N\in[L_1]:$ if $N'=\spp\zi_1\cup N$ and $\exists\, F\in\calf$
satisfying\\
 $\rule{2.5cm}{0cm}\min\{\langle\zi^{N'}_k,F\rangle:k=1,\ldots,n\}>\e$\\
$\rule{2.5cm}{0cm}$ and $F\cap\spp\zi^{N'}_1=F\cap\spp\zi_1=F_i\}$.

\noindent Each $A_i$ is an open set hence by infinite Ramsey theorem
we get that \tx  $i_0$ and $S_1\in[L_1]$ \st \fe $N\in[S_1]$ \txs
$F\in\calf$ with: 
\[\min\{\langle\zi^{N'}_k,F\rangle:k=1,\ldots,n\}>\e\]
 and
\[\spp\zi^{N'}_1\cap F=F_{i_0}.\]
Set $G_1=F_{i_0}$ and consider the set
\[\calf_{G_1}=\{F\in\calf:F\cap\spp\zi_1=G_1\}.\]
Then it is easy to see that $S_1,\z, \calf_{G_1}, n-1$ satisfy the
inductive assumptions hence \txs $N_1\in[S_1]$ \st
\[s(\calf_{G_1}[N_1])>(n-1)\cdot\om^\z.\]
As a  con\sq\ of this we get  that
$G_1\in\calf_{G_1}[N_1]^{(n-1)\cdot\om^\z}$.

Inductively choose $L\sps N_1\sps\cdots\sps N_k\sps\cdots$ and
$(G_k)_{k\in{\bf N}}$ \st

(i) $\max\spp\zi_k<\min N_k$

(ii) $G_k\sbs\spp\zi_k$, $G_k\in\calf$ and 
$\langle\zi_k,G_k\rangle >\e$

(iii) $G_k\in\calf_{G_k}[N_k]^{(n-1)\cdot\om^\z}$

\noindent where $\calf_{G_k}=\{F\in\calf:F\cap\spp\zi_k=G_k\}$.

The choice is done as in the case $\zi_1$.

Choose a set $N$ that is almost contained in $N_k$ for all $k\in\NN$.

\vspace{4pt}
{\sl Claim:} For every $k\in\NN$ the set 
$G_k$ belongs to $\calf_{G_k}[N]^{(n-1)\cdot\om^\z}$, 
where $\calf_{G_k}[N]$ is
defined as:
\[\calf_{G_k}[N]=\{F\in\calf[\spp\zi_k\cup N]:F\cap\spp\zi_k=G_k\}.\]

Indeed, set $N^k=\{m\in N:\max\spp\zi_k<m\}$.

Then since $N\stm N^k$ is finite, we get that
$G_k\in\calf_{G_k}[N]$ provided that $G_k\in\calf_{G_k}[N^k]$. 
Further,
$N^k$ is almost contained in $N_k$ and from (iii) and the fact that 
$(n-1)\cdot\om^\z$ is a limit ordinal we get that 
$G_k\in\calf_{G_k}[N]$.

To finish the proof of the lemma we prove the following.

\vspace{4pt}
{\sl Claim 2:} There exists $N'\in[N]$ \st $\calf[N']^{n\om^\z}\neq
\emptyset$.

Indeed, consider the family
\[\calg[N]=\{F:\exists\, k\in\nn\;\;{\rm with}\;\;
G_k\sbs N,\; F\sbs G_k\}.\]
It is easy to see that $\calg[N]$ is an adequate family. Further, for
$G_k\sbs N$ we have $\calf_{G_k}[N]\sbs\calf[N]$ hence \fe
$F\in\calg[N]$ we have that $F\in\calf[N]^{(n-1)\cdot\om^\xi}$.

So we get that $ \calg[N]\sbs\calf[N]^{(n-1)\cdot\om^\z}$.

Notice also that \fe $S\in[N]$ \txs $G_k\in\calg[N]$ \st
$\langle\zi^{S}_1,G_k\rangle >\e$.
Hence $\calg[N]$ is $(N,\z,\e)$ large and by the inductive assumption \txs
$N'\in[N]$ \st $\calg[N']^{\om^\z}\neq\emptyset$.
Hence 
\[\calf[N']^{n\cdot\om^\z}=\left[\calf[N']^{(n-1)\cdot\om^\z}
\right]^{\om^\z}\sps G[N']^{\om^\z}\neq\emptyset.\]

This completes the inductive proof of the lemma.

\sprh{\bf Completion of the proof of the proposition.}
Using the previous lemmas, for a given $L\in[M]$, we choose
$L\sps L_1\sps L_2\sps\cdots\sps L_n\sps\cdots$ \st
$s(\calf[L_n])>n\cdot\om^\z$. Then it is easy to see that if
$N$ is any set almost contained in $L_n$ for all $\nin$ then
$s(\calf[N])>n\cdot\om^\z$ for all $\nin$ and hence
$s(\calf[N])>\om^{\z+1}=\om^\xi$. 
The proof of the proposition is complete.

\sprh
We conclude this section with the following proposition.

\th{2.3.6 Proposition} Let $\xi<\om_1$, $M\in\N$, $\e>0$ and
$\calf$ be an adequate family. Suppose that \txs $L\in[M]$,
$L=(m_n)_{\nin}$ satisfying the property

for every $\nin$, $N\in[L]$, $n\leq\min N$

\centerline{$\sup_{F\in{\cal F}}\left\{\langle\xii^N_k,F\rangle:
k=1,2,\ldots,n\right\}>\e.$}

\noindent Then \txs $N\in[L]$ \st $s(\calf[N])>\om^{\xi+1}.$

\proof Notice that \calf\ satisfies the assumptions of the previous
lemma hence \txs a decreasing \sq\ $(L_n)_{\nin}$ of subsets of
$L$ \st $s(\calf[L_n])>n\cdot\om^\xi$.
Now, if $N$ is almost contained in $L_n$ for all $\nin$, it is easy
to see that $s(\calf[N])>\om^{\xi+1}$.

\sprt\sprh{\large\bf 2.4 The main results}

\sprt We pass now to give the statements and the proofs of the main results.

\th{2.4.1 Theorem} For $(x_n)_{\nin}$ weakly null \sq\ in a Banach space
$X$ and $\xi<\om_1$  exactly one of the following holds.

(a) For every $M\in[\nn]$ \txs $L\in[M]$ \st
\fe $P\in[L]$ the  \sq\ $(x_n)_{\nin}$ is $(P,\xi)$ summable

(b) There exists $M\in\N$ $M=(m_i)_{i\in\nn}$ \st
$(x_{m_i})_{i\in\nn}$ is an $\ell^1_{\xi+1}$ spreading model.

\sprh To prove the theorem we begin with the following lemma.

\th{2.4.2. Lemma} Assume that $F=(x_n)_{\nin}$ is a weakly null \sq\ and
$\xi<\om_1$ \st \fe $M\in\N$ \txs $N\in[M]$ with $(x_n)_{\nin}$
not $(N,\xi)$ summable. Then \txs $\e>0$ and $L\in\N$ \st
\fe 
$N\in[L]$ 
\[\overline{\lim}\,\left\| z^L_n\right\|>\e,\;\;{\rm 
where}\;\; z^L_n=\frac{\sum_{k=1}^n\xii^L_k\cdot F}{n}.\]

\proof We prove it for $M=\nn$.
The general case is similar.\\
For  given $\e>0$, $\nin$, we consider the set
\[ A_{\e,n}=\{M\in\N:\|z^M_k\|\leq\e\;\forall\, k\geq n\}.\]
Clearly each $A_{\e,n}$ is a closed set hence the set
$A_\e=\bigcupnf A_{\e,n}$ is a Ramsey set. Therefore \txs
$L_\e\in\N$ \st $[L_\e]\sbs A_\e$ or $[L_\e]\sbs\N\stm A_\e$.

If \txs some $\e>0$ and $L\in\N$ \st $[L]\sbs\N\stm A_\e$ then the
lemma has been proved.
Assume that this does not occur. Then inductively choose
$\nn\sps L_1\sps L_2\sps\cdots\sps L_n\sps\cdots$ \st
$[L_n]\sbs A_{1/n}$ and set $L$ any infinite set almost contained
in $L_n$ for all $\nin$. Set $L$ any set that is almost
contained in each $L_n$.

\vspace{8pt}
\noindent{\sl Claim 1:} If $N\in[L]$ then $(x_n)_{\nin}$ is $(N,\xi)$
summable.

Indeed, for any such $N$ and $\nin$ \txs $k_n\in\nn$ \st \fe
$k\in\nn$ with $k_n\leq k$ we have that 
 $\spp\xii^N_k\sbs L_n$. Therefore if
$N'=\bigcup_{k\geq k_n}\spp\xii^N_k$ then by the property P.4 we get
that $\xii^{N'}_s=\xii^N_{(s-1)+k_n}$.\\ 
Since $N'\in[L_n]$, \txs $s_0$ \st for all $s\geq s_0$
we have that $\|z^{N'}_s\|\leq\frac{1}{n}.$
But then \tx large $s_1$ \st \fe  $s>s_1$ we have
\beq
\left\|z^N_s\right\|\!\!\!=\!\!\!
\left\|\frac{\sum_{i=1}^s\xii^N_i\!\cdot\! F}{s}\right\|
&\!\!=\!\!&\left\|\frac{\sum_{i=1}^{k_n-1}\xii^N_i\!\cdot\! F}{s}+
\frac{\sum_{i=1}^s\xii^{N'}_i\!\cdot\! F}{s}-
\frac{\sum_{i=s+1-k_n}^s\xii^{N'}_i\!\cdot\! F}{s}\right\| \\
&\!\!<&\!\!\frac{k_n-1}{s}+\frac{1}{n}+\frac{k_n-1}{s}<\frac{2}{n}.
\eeq
This proves the Claim and it contradicts our assumptions.
Hence \txs $\e>0$ and $L\in\N$ \st $[L]$ is a subset of $\N\stm A_\e$
and this completes the proof of the lemma.

\th{2.4.3 Lemma} Let $F=(x_n)_{\nin}$ be  a
weakly null \sq. 
Suppose that for $\xi<\om_1$,  $M\in\N$ and $\e>0$ we have that
for all $N\in[M]$
\[\overline{\lim}\|z^N_n\|>\e\;\;{\rm for\; all}\;\; N\in[M],\;
{\rm where}\; z^N_n=\frac{\sumin\xii^N_i\cdot F}{n}.\]
Then \fe $L\in[M]$ we have:

(a) For every $\nin$ \txs $L_n\in[L]$ \st \fe $N\in[L_n]$
\[ a^N_n=\sup_{x^*\in B_{X^*}}\min\{x^*(\xii^N_k\cdot F):k=1,2,\ldots,n\}
>\frac{\e}{2}.\]

(b) There exists $N\in[L]$, $N=(m_n)_{\nin}$ such that:\\
\fe $N'\in[N]$ with $m_n\leq\min N'$ 
\[ a^{N'}_n>\frac{\e}{2}.\]

\proof (a) For a given $L\in[M]$ and $\nin$ we define a partition
of $[L]$ into $A_1,\; A_2$ such that:
\[ A_1=\left\{ n\in[L]:a_n^N>\frac{\e}{2}\right\},\;\;\;
A_2=[L]\stm A_1.\]
The set $A_1$ is a Ramsey set hence \txs $L_n$ \st either
$[L_n]\sbs A_1$ or $[L_n]\sbs A_2$. The first case proves
part (a) of the Lemma. We show that the second case does not occur.

Indeed, if $[L_n]\sbs A_2$ then we get that for 
$k_1<k_2<\cdots<k_n$ \txs $N\in[L_n]$ \st
$\xii^{L_n}_{k_1}=\xii^N_1,\ldots,\xii^{L_n}_{k_n}=\xii^N_n$
and since $N\in A_2,$ $a^N_n\leq\frac{\e}{2}$.
Choose $s$ large \st \txs $x^*\in B_{X^*}$ with
\[ x^*\left(\frac{\sum_{i=1}^s\xii^{L_n}_i\cdot F}{s}\right)>\e.\]
Then from the choice of $s$ we get that
\[ \#\left\{i:i\leq s,\: x^*(\xii^{L_k}_i\cdot F)>\frac{\e}{2}\right\}
\geq n.\]
 But then \txs $N\in[L_n]$ with $a^n_N>\frac{\e}{2}$,
a contradiction and the proof of part (a) is complete.

(b) Choose, inductively, a decreasing \sq\ $(L_n)_{\nin}$ \st
$L_n\in[L]$ and $L_n$ satisfies the requirement for the number $n$ 
of part (a). It is clear that any $N$ almost contained in $L_n$
for all $\nin$ is the desired set.

\sprh
Next we will prove two lemmas that will help us to reduce the proof
of the theorem to the case of the \sq\ $(\pi_n)_{n\in L}$ of the
natural coordinate projections of $\{0,1\}^{\nn}$ 
acting on an adequate family
\calf\ of finite subsets of $\nn$.

\th{2.4.4 Definition} Let $D$ be a weakly compact subset of
$c_0(\nn)$ and $\delta>0$. We set

$\calf_\de=\{ F\sbs\nn:\exists\, f\in D$ with $f(n)\geq\de\;\forall
\, n\in F\}$.

\th{2.4.5 Remark} The weak compactness of $D$ implies that
$\calf_\de$ is an adequate family of finite subsets of $\nn$.

\th{2.4.6 Notation} In the sequel we denote by

$\calf_\de[N]=\left\{ F:F\in\N^{<\om},\; \exists\, G\in\calf_\de:
G\cap N=F\right\}$.

\noindent Notice that  $\calf_\de[N]$ as a projection of
 $\calf_\de$ is also
a compact family.

\sprh The next Lemma is a con\sq\ of Lemma 1.2.

\th{2.4.7 Lemma} Let $D$ be a weakly compact subset of
$c_0(\nn)$. Then \fe $\de>0$, $\e>0$ and $M\in\N$ \txs $N\in[M]$
\st \fe $F\in\calf_\de[N]$ and $G\sbs F$ \txs $f\in D$ \st

(i) $\min\{f(n):n\in G\}\geq(1-\e)\de$

(ii) $\sum_{n\not\in G}|f(n)|<\e\cdot\de$.

\proof From Lemma 1.2 \fe $M\in\N$ \txs $N\in[M]$ \st \fe $k\in\nn$,
$\max F\leq k$ \txs $f_k\in D$

(i) $\min\{f_k(n):n\in G\}>(1-\e)\de$

(ii) $\sum_{n=1 \atop n\not\in G}^k|f_k(n)|<\e\cdot\de$.

\noindent The desired $f$ is weak limit of any weakly convergent sub\sq\ of
$(f_k)_{k\in{\bf N}}$.

\sprh{\bf 2.4.8 Lemma} {\sl (Reduction lemma)} Let $H=(x_n)_{\nin}$ be a
weakly null \sq\ in a Banach space with $\|x_n\|\leq 1$.
Then \fe $\de>0$ and $\e>0$ \txs an adequate family \calf\ of finite subsets
of \nn\ and a function $f:B_{X^*}\lra\calf$ such that:

For every $M\in\N$ \txs $N\in\N$ satisfying the following properties

(a) If $A\in S^+_{\ell^1}$, $\spp A\sbs N$ then for every 
$x^*\in B_{X^*}\langle x^*,A\cdot H\rangle >\de$ we have
$\langle A,f(x^*)\rangle >\frac{\de^2}{4}$.

(b) If $A\in S_{\ell^1}$, $\spp A\sbs N$ and $F\in\calf$ \st
$\langle A,F\rangle\geq\frac{1}{2}$ then
$\|A\cdot H\|\geq\frac{\e\cdot\de}{4}$.

\proof We start by noticing that if
$A\in S^+_{\ell^1}$ and $x^*\in B_{X^*}$ \st 
$A=(a_n)_{\nin}$ and $x^*(A\cdot H)>\de$ then for 
$F=\left\{\nin:x^*(x_n)>\frac{\de}{2}\right\}$ we get that
$\sum_{n\in F}a_n>\frac{\de}{2}$. Hence 
$\langle A,F\rangle >\frac{\de^2}{4}$.

Since $(x_n)_{\nin}$ is a weakly null \sq , the set
\[ D=\{(x^*(x_n))_{\nin}:x^*\in B_{X^*}\}\]
 is a weakly compact subset
of $c_0(\nn)$. Applying  Lemma 2.4.7 for $\frac{\de}{2}$,
$\frac{\e}{2}$ and $M\in\N$ we find $N\in[M]$ satisfying
properties (i) and (ii) of that lemma. 
We let \calf\ be the adequate
family defined as $\calf_{\frac{\de}{2}}$.
We also define  $f:B_{X^*}\lra\calf$ by the rule
$f(x^*)=\left\{\nin:x^*(x_n)>\frac{\de}{2}\right\}.$
Using our note at the beginning of the proof, we get that property (a)
holds \fe $N\in[M]$.
To see property (b), suppose that $A\in S_{\ell^1}$ with
$\spp A\sbs N$ and $F\in\calf$ \st 
$\langle A,F\rangle \geq\e$. Then we may assume that
$F\sbs N\cap\{\nin:a_n>0\}$ and by the definition of \calf\
\txs $G\in\calf_{\frac{\de}{2}}$ \st $F\sbs G\cap N$.
Then \txs $x^*\in B_{X^*}$ \st

(i) $\min\{x^*(x_n):n\in F\}\geq\e\lf(1-\f{\e}{2}\r)\f{\de}{2}$

(ii) $\sum_{n\not\in F}|x^*(x_n)|<\frac{\e}{2}$.

\noindent From (i), (ii) and the fact that 
$\langle A,F\rangle \geq\e$ we get that 
$\|A\cdot H\|>\e\cdot\frac{\de}{4}.$
The proof is complete.

\sprh{\bf Proof of the theorem:}
We prove first that the negation of (a) implies (b).
Suppose that $\|x_n\|\leq 1$ and for a given $\xi<\om_1$
the case (a) does not occur.
Then from Lemma 2.4.2 \txs $M\in\N$ and $\de>0$ \st
$\overline{\lim}\|z^L_n\|>2\de$ for all $L\in[M]$.
Going to a subset  of $M$ if it necessary, we may assume
that part (b) of Lemma 2.4.3 is also satisfied for $M$.

Consider the family \calf\ defined in Lemma 2.4.8 for the \sq\ 
$(x_n)_{\nin}$ and the number $\de$.
Let  $N\in[M]$ \st (a) and (b) in Lemma 2.4.8 are
satisfied. Property (a) in connection with the fact that $N$
satisfies the conclusion of Lemma 2.4.3 show that the assumptions of
Proposition 2.3.6 are fulfilled.
 Hence \txs $N'\in[N]$
\st $s(\calf[N'])>\om^{\xi+1}$.
>From Theorem 2.2.6 there exists $N_1\in[N']$ \st
$N_1=(m_i)_{i\in {\bf N}}$ and \fe $F\in\calf_{\xi+1}$ the set
$\{m_i:i\in F\}\in\calf$.

\vspace{3pt}{\sl Claim:} For every $F\in\calf_{\xi+1}$
$\|\sum_{i\in F}a_ix_{m_i}\|\geq\frac{\de}{8}\sum_{i\in F}|a_i|$.

Indeed, by standard arguments, it is enough to show it for 
$(a_i)_{i\in F}\in S_{\ell^1}$.
If $(a_i)_{i\in F}\in S_{\ell^1}$ then either
$\sum\{a_i:i\in F,\: a_i>0\}\geq\frac{1}{2}$
or $\sum\{a_i:i\in F,\: a_i<0\}\leq -\frac{1}{2}$.
We assume that the first case occurs. Otherwise we consider the
$(b_i)_{i\in F}$ \st $b_i=-a_i$ for all $i\in F$.
Set $F'=\{i\in F:a_i>0\}$; then clearly
$\langle A,F'\rangle\geq\frac{1}{2}$ and hence
$\|A\cdot(x_n)_{\nin}\|>\frac{\de}{8}$,
which proves the claim.
The proof is complete.

\sprh We pass now to show that parts (a) and (b) of Theorem 2.4.1
are mutually exclusive.

\th{2.4.9 Proposition} Let $(x_n)_{\nin}$ be a weakly null \sq\ in a
Banach space $X$. If $\xi<\om_1$, $M\in\N$ and $\de>0$ are such that
$M=(m_i)_{i\in\nn}$ and 
\[\|\sum_{i\in F}a_ix_{m_i}\|\geq\de\cdot\sum_{i\in F}
|a_i|\;\; {\rm for\; every}\;\; F\in\calf_{\xi+1}\]
 then \txs $L\in[M]$ \st
\fe $P\in[L]$ $(x_n)_{\nin}$ is not $(P,\xi)$ summable.

\proof Consider the adequate family \calf\ defined in Reduction
lemma (Lemma 2.4.8) for the \sq\ $(x_n)_{\nin}$
the number $\de$ in our assumptions and
$\e=\de_{\xi+1}$ (Proposition 2.1.10).
Find  $N\in[M]$ \st conditions (a), (b) of Reduction Lemma
are fulfilled. Denote by 
$\lf(\lf[\xii^M\r]^P_n\r)_{\nin}$, $P\in[N]$,
the summabilidy methods defined by the rule
$\lf[\xii^M\r]^P_n=(a_{m_i})_{i\in\nn}$ where $\xii^P_n=(a_i)_{i\in\nn}$
and  $a_{m_i}=a_i$.
 Then \calf\ is
$\lf( N,\xii+1^M,\f{\de^2}{4}\r)$ large hence \txs
$L'\in\N$ \st $\cf^{L'}_{\xi+1}$ is a subfamily of \calf\ and
hence by Proposition 2.1.11 \txs $L''\in[L']$ \st \fe $P\in[L'']$,
$\nin$ \txs $G\in\cf^{L'}_{\xi+1}$ \st 
$\lf\langle(\mbox{\boldmath$\xi+1$})^P_n,G\r\rangle >\de_{\xi+1}>0$.
Choose, as in Lemma 2.4.3 (part (b)), an $L\in[L'']$ \st
$L=(l_n)_{\nin}$ and \fe $P\in[L]$
$n\leq k_1<k_2<\cdots<k_n$ \txs $G\in\cf$ \st
$\lf\langle\xii^P_{k_i},G\r\rangle >\f{\de_{\xi+1}}{2}$.
Then by part (b) of Reduction Lemma \txs $x^*\in B_{X^*}$ \st
$\lf\langle\xii^P_{k_i}\cdot H,X^*\r\rangle >\f{\de\cdot\de_{\xi+1}}{8}$
 where $H=\langle x_n\rangle_{\nin}$ and
$(k_i)_{i=1}^n,\; P$ as above.
It is clear now that \fe $P\in[L]$ the \sq\ $(x_n)_{\nin}$ is not
$(P,\xi)$ summable.

\th{2.4.10 Remark} The above Proposition immediately shows that parts
(a) and (b) in the Theorem 2.4.1 are mutually exclusive.

\sprh For the sequel we need the following result proved in [A-A].

\th{2.4.11 Proposition} Let $X$ be a Banach space and $(x_n)_{n\in\nn}$ a weakly
null \sq\ in $X$.

(a) There exists $\xi<\om_1$ \st for all $\z<\om_1$,
$\xi\leq\z$ $(x_n)_{n\in\nn}$ does not contain  a sub\sq\
which is $\ell^1_\z$ spreading model.

(b) If $\ell^1$ does not embed into $X$ then \txs $\xi<\om_1$ \st
\fe $\z<\om_1$, $\xi\leq\z$ and any bounded \sq\ $(x_n)_{n\in\nn}$
there is no sub\sq\ of $(x_n)_{n\in\nn}$ which is $\ell^1_\z$
spreading model.

{\em Sketch of proof:} The proof of this Proposition part (a) follows from
the fact that  
\[\calt_\e=\lf\{ (x_{n_1},\ldots, x_{n_k}):
\lf\|\sum_{i=1}^ka_ix_{n_i}\r\|\geq\e\cdot\sum_{i=1}^k|a_i|\r\}\]
ordered in the usual manner is a well-founded tree.
If not, the \sq\ $(x_n)_{n\in\nn}$ should contain a sub\sq\ equivalent to $\ell^1$
basis that contradicts the weak nullness of $(x_n)_{n\in\nn}$.
Therefore the height of $\calt_\e$, denoted by $o(\calt_\e)$, is a
countable ordinal $\xi_\e$. Further, if $(x_n)_{n\in\nn}$ has
a sub\sq\ that is $\ell^1_\xi$ spreading model with constant
$\de_\xi>\e$ then $\om^\xi\leq\xi_\e$.
So if $\xi_0=\sup\{\xi_\e, \e>0\}$ then every $\xi<\om_1$ \st
$(x_n)_{n\in\nn}$ has a sub\sq\ which is $\ell^1_\xi$ spreading model
should satisfy $\om^\xi\leq\xi_0$ and this proves the result
for part (a).

The proof of part (b) is the same and uses the technique developed
by Bourgain [B].

\sprt{\bf The Banach - Saks index}

{\bf 2.4.12 Definition} Let $X$ be a Banach space and $(x_n)_{n\in\nn}$ a weakly null
\sq\ in $X$.

(a) The Banach-Saks index of $(x_n)_{n\in\nn}$ denoted by $
BS[(x_n)_{n\in\nn}]$ is the
least ordinal $\xi$ \st there is no sub\sq\ of $(x_n)_{n\in\nn}$ which
is $\ell^1_\xi$ spreading model.

(b) If $X$ is a Banach space not containing $\ell^1(\nn)$ then we
denote by $BS[X]$ the least ordinal $\xi$ \st no bounded \sq\
$(x_n)_{n\in\nn}$ in $X$ is an $\ell^1_\xi$ spreading model.

\th{2.4.13 Theorem} Let $H=(x_n)_{\nin}$ be a weakly null \sq\ with\\
$BS[(x_n)_{\nin}]=\xi$.
Then: $\xi$ is the unique ordinal satisfying the following

(a) For every $M\in[\nn]$ \txs $L\in[M]$ \st \fe
$P\in[L]$ $\lim_{\nin}\|\xii^P_n\cdot H\|=0$.

(b) For every $\z<\xi$ \txs $L_\z\in\N$ \st $L_\z=(n_i)_{i\in\nn}$ and
$(x_{n_i})_{i\in\nn}$ is an $\ell^1_\z$ spreading model.

(c) If $\xi=\z+1$ \txs $\e>0$ and $L\in\N$ \st for all
$P\in[L]$ $\lf\|\zi^P_n\cdot H\r\|>\e$ and
$(\lf(\zi^P_n\cdot H\r)_{\nin}$ is Cesaro summable.

\proof (a) For $L\in\N$ and $\nin$ we define a partition of $[L]$
into sets $A$, $B$ by the rule
$A=\{ P:\|\xii^P_1\cdot H\|\leq\e\}$ and $B=[L]\stm A$.
It is easy to see that $A$ is a closed subset of $[L]$ hence
 by infinite Ramsey Theorem \txs $N\in[L]$ \st either 
$[N]\sbs A$ or $[N]\sbs B$.
If the second case holds then by Reduction lemmma we get that
$(x_n)_{n\in\nn}$ has a sub\sq\ which is $\ell^1_\xi$ spreading model,
a contradiction. Hence $[N]\sbs A$.
Choose, inductively, $M\sps L_1\sps L_2\sps\cdots\sps L_n\sps\cdots$
\st \fe $P\in[L_n]$
$\|\xii^P_1\cdot H\|<\frac{1}{n}$ and set
$L=(l_n)_{n\in\nn}$ \st $l_n\in L_n$. Then it is easy to see $L$
satisfies the conclusion of the first part of the theorem.

(b) It follows from the definition of $\xi$.

(c) Suppose now that $\xi=\z+1$. Then, by the definition of
$BS[(x_n)_{n\in\nn}]$, \txs $M\in\N$ \st $(x_n)_{n\in M}$ is an
$\ell^1_\z$ spreading model. Then by part (a) of Theorem 2.4.1
there exists $N\in[M]$ \st \fe $P\in[N]$
$(x_n)_{n\in\nn}$ is $(P,\xi)$ summable and finally from
Proposition 2.4.9 \txs $L\in[N]$ and $\e>0$ \st \fe $P\in[L]$, 
$\nin$, $\|\zi^P_n\cdot H\|\geq\e$. This proves part (c) and the
proof is complete.

\th{2.4.14 Remark} (i) The first part of the above Theorem is satisfied by
any normalized \sq\ in Tsirelson's space.
Any such \sq\ has Banach-Saks index equal to $\om$.

(ii) The third part gives a complete answer in the following question
posed by the first named author:
For what weakly null \sq s there exists a \sq\ $(y_n)_{n\in\nn}$ of
block convex combinations \st $\|y_n\|>\e$ and
$(y_n)_{n\in\nn}$ is Cesaro summable. 

\sprh We conclude this Section with the following corollaries.
Their proofs follow easily from the previous theorems.

\th{2.4.15 Corollary} For every separable reflexive Banach space $X$
\txs a unique ordinal $\xi<\om_1$ \st

(i) For all ordinals $\z\geq \xi$ the space $X$ has $\z$-BS.

(ii) For every $\z<\xi$ the space $X$ fails $\z$-BS.

\th{2.4.16 Corollary} If $X$ is a separable Banach space not containing
isomorphically $\ell^1$ then \tx a unique ordinal 
$\xi<\om_1$ \st

(i) For all ordinals $\z\geq \xi$ the space $X$ has w $\z$-BS.

(ii) For every $\z<\xi$ the space $X$ fails w $\z$-BS.

\sprt\sprt
\[{\bf REFERENCES }\]
\small
\begin{description}
\parsep =0pt
\itemsep =0pt
\parskip=0pt

\bb{[Al--Ar]}] {\sc D. Alspach} and {S. Argyros,}
{\em Complexity of Weakly null \sq s,}
Dissertations Mathematicae, CCCXXI, pp. 1-44, 1992.
\bb{[A--O]}] {\sc D. Alspach, E. Odell,}
{\em Averaging Weakly Null sequences,}
Lecture Notes in Math. 1332, Springer, Berlin 1988.
\bb{[A--D]}]
  {\sc S. Argyros} and {\sc I. Deliyanni,}  
{\em Examples of asymptotic $\ell^1$ Banach spaces} (preprint)
\bb{[B]}] {\sc Z. Bourgain,} {\em Convergent \sq s of continuous
functions,} Bull. Soc. Math. Belg. Ser. B 32 (1980), 235-249.
\bb{[Ell]}] {\sc E. Ellenduck,} {\em A new proof that analytic sets 
are Ramsey,} J. Symbolic Logic 39(1974), 163-165.
\bb{[E]}] {\sc J. Elton,} Thesis, Yale University
\bb{[E--M]}] {\sc P. Erd\"os} and {\sc M. Magidor,}
{\em A note on  regular methods of summability and the Banach-Saks
property}, Proc. Amer. Math. Soc. 59(1976), 232-234.
\bb{[G--P]}] {\sc F. Galvin -- K. Prikry,}
{\em Borel sets and Ramsey's theorem,}
J. Symbolic Logic 38 (1973) 193--198.
\bb{[G--M]}] {\sc W. T. Gowers} and {\sc B. Maurey,}
{\em The unconditional basic \sq\ problem,}
Journal of AMS 6, (1993, 851-874.
\bb{[M--R]}] {\sc B. Maurey} and {\sc H. Rosenthal,}
{\em Normalized weakly null \sq s with no unconditional sub\sq,}
Studia Math. 61 (1977) 77-98.
\bb{[M--N]}] {\sc S. Mercourakis, S. Negrepontis,} {\em Banach Spaces and
Topology II,} Recent Progress in General Topology, M. Husek and
J. Vaan Mill (eds.), Elsevier Sciences Publishers, 1992.
\bb{[M]}] {\sc S. Mercourakis,} {\em On Cesaro summable \sq s of continuous
functions,} Mathematica, 42(1995), 87-104.
\bb{[N--W]}] {\sc C. St. J. A. Nash-Williams,} 
{\em On well quasi-ordering transfinite \sq s,}
Proc. Camb. Phil. Soc. 61 (1965), 33-39.
\bb{[O$_1$]}] {\sc E. Odell,}
{\em Applications of Ramsey theorems to Banach space theory},
Notes in Banach Spaces (H. E. Lacey, ed.), Univ. of Texas Press,
1980, 379-404.
\bb{[O$_2$]}] {\sc E. Odell,} {\em On Schreier unconditional \sq s,}
Contemporary Math.\\ 144(1993), 197-201.
\bb{[Sch]}] {\sc J. Schreier,} {\em Ein Gegenbeispiel zur Theorie
der schwachen Konvergenz,} Studia Math. 2(1930), 58-62.
\bb{[Si]}] {\sc J. Silver,} {\em Every analytic set is Ramsey,}
J. Symbolic Logic 35(1970), 60-64.
\bb{[T]}] {\sc B. S. Tsirelson,} {\em  Not every Banach space contains
$\ell_p$ or $c_0$}, Funct. Anal. Appl. 8, 1974,  138-141.
\end{description}

\sprt
\mbox{}\hfill\begin{minipage}[t]{6cm}
Department of Mathematics\\
University of Athens\\
Athens 15784, Greece\\
{\tt e-mail: sargyros @ atlas.uoa.}\\
\mbox{}\hfill {\tt ariadne - t.gr$\;\;\;\;\;\;\;$}
\end{minipage}

\end{document}